\input amstex
\documentstyle{amsppt}
\input bull-ppt
\keyedby{bull284e/lic}
\define\Aut{\mathop{\roman {Aut}}\nolimits}
\define\CC{{\bold C}}
\define\chr{\mathop{\roman {char}}\nolimits}
\define\coker{\mathop{\roman {coker}}\nolimits}
\define\Cl{\mathop{\roman {Cl}}\nolimits}
\define\Dd{{\scr D}}
\define\den{\mathop{\roman {den}}\nolimits}
\define\End{\mathop{\roman {End}}\nolimits}
\define\FF{{\bold F}}
\define\Gal{\mathop{\roman {Gal}}\nolimits}

\define\Hom{\mathop{\roman{ Hom}}\nolimits}
\define\Oo{{\scr O}}
\define\ord{\mathop{\roman{ ord}}\nolimits}

\define\reg{\mathop{\roman {reg}}\nolimits}
\define\QQ{{\bold Q}}
\define\RR{{\bold R}}
\define\Spec{\mathop{\roman{ Spec}}\nolimits}
\define\Tr{\mathop{\roman{ Tr}}\nolimits}
\define\vl{\mathop{\roman{ vol}}\nolimits}
\define\ZZ{{\bold Z}}

\define\xaa{{\germ a}}
\define\bbx{{\germ b}}
\define\ccx{{\germ c}}

\define\NNx{{\germ N}}
\define\ppx{{\germ p}}
\define\qqx{{\germ q}}
\define\rrx{{\germ r}}
\define\yyx{{\germ y}}

\define\ADHU{1}  
\define\ADLE{2}  
\define\ADPR{3}  
\define\ARCH{4}  
\define\ATMO{5}  
\define\BACH{6}  
\define\BASH{7}  
\define\BERW{8}  
\define\BOSH{9}  
\define\BOHU{10} 
\define\BRKN{11} 
\define\BUCA{12} 
\define\BUCB{13} 
\define\BULE{14} 
\define\BUSH{15} 
\define\BUWI{16} 
\define\BULP{17} 
\define\CAME{18} 
\define\CASS{19} 
\define\CAFR{20} 
\define\CHIA{21} 
\define\CHIB{22} 
\define\COHE{23} 
\define\COLA{24} 
\define\COLH{25} 
\define\FOMK{26} 
\define\GORD{27} 
\define\HMCA{28} 
\define\HMCB{29} 
\define\KANN{30} 
\define\KANT{31} 
\define\KALU{32} 
\define\KNUT{33} 
\define\LANB{34} 
\define\LANA{35} 
\define\LAMI{36} 
\define\LANG{37} 
\define\AKLT{38} 
\define\AKFA{39} 
\define\AKFB{40} 
\define\LLAN{41} 
\define\LLLF{42} 
\define\LLMF{43} 
\define\LLMP{44} 
\define\LERE{45} 
\define\LEGA{46} 
\define\LEFL{47} 
\define\LEFI{48} 
\define\LEPO{49} 
\define\LETI{50} 
\define\LOVO{51} 
\define\MINE{52} 
\define\ODDL{53} 
\define\PALF{54} 
\define\POHS{55} 
\define\POZA{56} 
\define\POFR{57} 
\define\QUEM{58} 
\define\SAND{59} 
\define\SCHI{60} 
\define\SCLT{61} 
\define\SCEC{62} 
\define\SEAA{63} 
\define\SEMW{64} 
\define\SHAN{65} 
\define\SHOA{66} 
\define\SHOB{67} 
\define\SIEG{68} 
\define\STAU{69} 
\define\SZPI{70} 
\define\TEIT{71} 
\define\TZDW{72} 
\define\VDLI{73} 
\define\VEMB{74} 
\define\WEIS{75} 
\define\ZANT{76} 
\define\ZASS{77} 
\define\ZIMM{78} 

\topmatter
\cvol{26}
\cvolyear{1992}
\cmonth{April}
\cyear{1992}
\cvolno{2}
\cpgs{211-244}
\title Algorithms in Algebraic Number Theory \endtitle
\author H. W. Lenstra, Jr.\endauthor
\shortauthor{H. W. Lenstra, Jr.}
\shorttitle{Algorithms in algebraic Number Theory}
\address Department of Mathematics, University of 
California,
Berkeley, California 94720 \endaddress
\date October 11, 1991\enddate
\subjclassrev{Primary 11Y16, 11Y40}
\thanks The author was supported by NSF under Grant No.\ 
DMS 90-02939\endthanks
\keywords Algebraic number theory, algorithms, complexity 
theory\endkeywords
\abstract In this paper we discuss the basic problems of 
algorithmic algebraic
number theory. The emphasis is on aspects that are of 
interest from a
purely mathematical point of view, and practical issues 
are largely
disregarded. We describe what has been done and, more 
importantly,
what remains to be done in the area. We hope to show that 
the study
of algorithms not only increases our understanding of 
algebraic number
fields but also stimulates our curiosity about them. The 
discussion is
concentrated of three topics: the determination of Galois 
groups, the
determination of the ring of integers of an algebraic 
number field,
and the computation of the group of units and the class 
group of that
ring of integers. \endabstract

\endtopmatter

\document



\heading1. Introduction\endheading
The main interest of algorithms in algebraic number theory 
is that they
provide number theorists with a means of satisfying their 
professional
curiosity. The praise of numerical experimentation in 
number theoretic
research is as widely sung as purely numerological 
investigations are
indulged in, and for both activities good algorithms are 
indispensable.
What makes an algorithm {\it good\/} unfortunately defies 
definition---too
many extra-mathematical factors affect its practical 
performance, such as
the skill of the person responsible for its execution and 
the
characteristics of the machine that may be used.
\fn""{This paper was given as a Progress in Mathematics 
Lecture at
the August 8--10, 1991 meeting of the American 
Mathematical Society
in Orono, Maine.}

The present paper addresses itself not to the researcher 
who is looking for
a collection of well-tested computational methods for use 
on his recently
acquired personal computer. Rather, the intended reader is 
the perhaps
imaginary pure mathematician who feels that he makes the 
most of his talents
by staying away from computing equipment. It will be 
argued that even from
this perspective the study of algorithms, when considered 
as objects of
research rather than as tools, offers rich rewards of a 
theoretical nature.

The problems in pure mathematics that arise in connection 
with algorithms
have all the virtues of good problems. They are of such a 
distinctly
fundamental nature that one is often surprised to discover 
that they have
not been considered earlier, which happens even in 
well-trodden areas of
mathematics; and even in areas that are believed to be 
well-understood it
occurs frequently that the existing theory offers no ready 
solutions,
fundamental though the problems may be. Solutions that 
have been found often
need tools that at first sight seem foreign to the 
statement of the problem.

Algebraic number theory has in recent times been applied 
to the solution of
algorithmic problems that, in their formulations, do not 
refer to algebraic
number theory at all. That this occurs in the context of 
solving diophantine
equations (see, e.g., [\TZDW]) does not come as a 
surprise, since these
lie at the very roots of algebraic number theory. A better 
example is
furnished by the seemingly elementary problem of 
decomposing integers into
prime factors. Among the ingredients that make modern 
primality tests work one
may mention reciprocity laws in cyclotomic fields (see 
[\ADPR, \COLH, \COLA]),
arithmetic in cyclic fields (see [\LEGA, \BOHU]), the 
construction of
Hilbert class fields of imaginary quadratic fields 
[\ATMO], and class number
estimates of fourth degree CM-fields [\ADHU]. The best 
rigorously proved time
bound for integer factorization is achieved by an 
algorithm that depends on
quadratic fields (see [\LEPO]), and the currently most 
promising practical
approach to the same problem, the {\it number field 
sieve\/} (see [\BULP,
\LLMF, \LLMP]), employs ``random'' number fields of which 
the discriminants
are so huge that many traditional computational methods 
become totally
inapplicable.  The analysis of many algorithms related to 
algebraic number
fields seriously challenges our theoretical understanding, 
and one is often
forced to argue on the basis of heuristic assumptions that 
are formulated for
the occasion. It is considered a relief when one runs into 
a standard
conjecture such as the generalized Riemann hypothesis (as 
in [\BACH, \BUSH]) or
Leopoldt's conjecture on the nonvanishing of the $p$-adic 
regulator [\SCHI].

In this paper we will consider algorithms in algebraic 
number theory for
their own sake rather than with a view to any of the above 
applications.
The discussion will be concentrated on three basic 
algorithmic questions that
one may ask about algebraic number fields, namely, how to 
determine the
Galois group of the normal closure of the field, or, more 
generally, of any
polynomial over any algebraic number field; how to find 
the ring of integers
of the field; and how to determine the unit group and the 
ideal class group
of that ring of integers. These are precisely the subjects 
that are discussed
in {\it Algorithmic algebraic number theory\/} by M. Pohst 
and H. Zassenhaus
(Cambridge, 1989), but our point of view is completely 
different. Pohst and
Zassenhaus present algorithms that ``yield good to 
excellent results for
number fields of small degree and not too large 
discriminant'' [\POZA,
Preface], but our attitude will be decidedly and 
exclusively asymptotic. For
the purposes
of the present paper one algorithm is considered better 
than another if,
for each positive real number $N$, it is at least $N$ 
times as fast for all
but finitely many values of the input data. It is clear 
that with this
attitude we can make no claims concerning the practical 
applicability of any
of the results that are achieved. In fact, following 
Archimedes [\ARCH] one
should be able, on the basis of current physical 
knowledge, to find an
upper estimate for all sets of numerical input data to 
which any algorithm
will ever be applied, and an algorithm that is faster in 
all those finitely
many instances may still be worse in our sense. 

To some people the above attitude may seem absurd. To the 
intended reader,
who is never going to apply any algorithm anyway, it comes 
as a liberation
and a relief. Once he explicitly gives up all practical 
claims he will
realize that he can occupy himself with algorithms without 
having to fear
the bad dreams caused by the messy details and dirty 
tricks that stand
between an elegant algorithmic idea and its practical 
implementation. He
will find himself in the platonic paradise of pure 
mathematics, where a
conceptual and concise version of an algorithm is valued 
more highly than
an {ad hoc\/} device that speeds it up by a factor of ten 
and where
words have precise meanings that do not change with the 
changing world. He
will never need to enter the dark factories that in his 
imagination are
populated by applied mathematicians, where boxes full of 
numbers that they
call matrices are carried around and where true electronic 
computers are
fed with proliferating triple indices. And in his 
innermost self he will know
that in the end his own work will turn out to have the 
widest application
range, exactly because it was not done with any specific 
application in mind.

There is a small price to be paid for admission to this 
paradise. Algorithms
and their running times can only be investigated 
mathematically if they are
given exact definitions, and this can apparently be done 
only if one employs
the terminology of {\it theoretical computer science}, 
which our intended
reader unfortunately does not feel comfortable with 
either. It is only out of
respect for his feelings that I have not called this paper 
{\it Complexity of
algorithms in algebraic number theory}, which would have 
described its
contents more accurately.

Although it is, from a rigorous mathematical point of 
view, desirable that
I define what I mean by an algorithm and its running time, 
I will not do so.
My main excuse is that I do not know these definitions 
myself. Even worse, I
never saw a treatment of the appropriate theory that is 
precise, elegant, and
convenient to work with. It would be a laudable enterprise 
to fill this
apparent gap in the literature. In
the meantime, I am happy to show by example that one can 
avoid paying the
admission price, just as not all algebraists are experts 
on set theory or
algebraic geometers on category theory. The intuitive 
understanding that one
has of algorithms and running times, or of sets and 
categories, is amply
sufficient. Exact definitions appear to be necessary only 
when one wishes to
prove that algorithms with certain properties do {\it 
not\/} exist, and
theoretical computer science is notoriously lacking in 
such negative results.
The reader who wishes to provide his own definitions may 
wish to consult
[\VEMB] for an account of the pitfalls to be avoided. He 
should bear in mind
that all theorems in the present paper should become 
formal consequences of
his definitions, which makes his task particularly academic.

My intended reader may have another allergy, namely, for 
{\it constructive
mathematics}, in which purely existential proofs and the 
law of the excluded
middle are not accepted. This has only a superficial 
relationship to
algorithmic mathematics. Of course, it often happens that 
one can obtain a
good algorithm by just transcribing an essentially 
constructive proof, but
such algorithms do not tend to be the most interesting 
ones; many of them
are mentioned in \S 2. In the design and analysis of 
algorithms one
gladly invokes all the help that existing pure mathematics 
has to offer and
often some not-yet-existing mathematics as well.

For an account of algorithms in algebraic number theory 
that emphasizes the
practical aspects rather than complexity issues we refer 
to the
forthcoming book by Cohen [\COHE].

In \S 2 we cover the basic terminology and the basic 
auxiliary
results to be used in later sections. In particular, we 
discuss several
fundamental questions that, unlike integer factorization, 
admit a
satisfactory algorithmic treatment. These include 
questions related to
finitely generated abelian groups, to finite fields, and 
to the factorization
of polynomials over number fields.

Section 3 is devoted to the problem of determining Galois 
groups. We review
the little that has been done on the complexity of this 
problem, including the
pretty result of Landau and Miller [\LAMI] that 
solvability by radicals can
be decided efficiently. We also point out several 
directions for further
research.

In \S4 we discuss the problem of determining the ring of 
integers of a
given algebraic number field. The main result is a 
negative one---the problem
is in many ways equivalent to the problem of finding the 
largest square
factor of a given positive integer, which is intractable 
at present.
Nevertheless, we will see that one can get quite close. 
There is an
interesting connection with the resolution of plane curve 
singularities that
remains to be exploited.

Section 5 considers the problem of determining the unit 
group $\Oo^*$ and the
ideal class group $\Cl\Oo$ of the ring of integers $\Oo$ 
of a given number
field. Showing that these are effectively computable is 
not entirely trivial,
and since most textbooks are silent on this point, I treat 
it in some detail.
We shall see that the existing complexity estimates for 
this problem still
leave room for improvement, and what we have to say is far 
from conclusive.
In \S6 we prove a few explicit bounds
concerning units and class groups that are needed in \S 5. 
Several
results in these two sections could have been formulated 
in terms of the
divisor class group $\operatorname{Pic}_c\Oo$ that  
appears in Arakelov theory (see
[\SZPI, \S I]) and that already appeared in the context of
algorithms (see [\SHAN, \LERE]). Knowing the group 
$\operatorname{Pic}_c\Oo$ is equivalent
to knowing
both $\Oo^*$ and $\Cl\Oo$, which may explain why 
algorithms for computing
$\Oo^*$ and algorithms for computing $\Cl\Oo$ are often 
inextricably linked.
It also explains why, contrary to many authors in the 
field, I prefer to
think of determining $\Oo^*$ and determining $\Cl\Oo$ as a 
single problem.

The three basic questions that are addressed in this {\it 
progress report\/}
still offer ample opportunities for additional progress. 
Among the many other
algorithmic questions in algebraic number theory that 
merit attention we
mention the problem of tabulating number fields, problems 
from class field
theory such as the calculation of Artin symbols, problems 
concerning quadratic
forms, and the analogues of all problems discussed in this 
paper for function
fields of curves over finite fields.
\heading 2. Preliminaries\endheading

\subheading{{\rm2.1.}  Algorithms and complexity}
It is assumed that the reader has an intuitive 
understanding of the
notion of an {\it algorithm\/} as being a recipe that 
given one finite
sequence of nonnegative integers, called the {\it input\/} 
data, produces
another, called the {\it output}. Formally, an algorithm 
may be defined as a
{\it Turing machine}, but for several of our results it is 
better to
choose as our ``machine model'' an idealized computer that 
is more realistic
with respect to its {\it running time}, which is another 
intuitively clear
notion that we do not define. We refer to [\VEMB] and the 
literature given
there for a further discussion of these points.

The {\it length\/} of a finite sequence of nonnegative 
integers $n_1$, $n_2$,
\dots, $n_t$ is defined to be $\sum_{i=1}^t\log (n_i+2)$. 
It must informally
be thought of as proportional to the number of bits needed 
to spell out the
$n_i$ in binary. By analyzing the {\it complexity\/} of an 
algorithm we mean
in this paper finding a reasonably sharp upper bound for 
the running time of
the algorithm expressed as a function of the length of the 
input data. This
should, more precisely, be called {\it time\/} complexity, 
to distinguish it
from {\it space\/} complexity. An algorithm is said to be 
{\it
polynomial-time\/} or {\it good\/} if its running time is 
$(l+2)^{O(1)}$, where
$l$ is the length of the input. Studying the complexity of 
a {\it problem\/}
means finding an algorithm for that problem of the 
smallest possible
complexity. In the present paper we consider the 
complexity analysis complete
when a good algorithm for a problem has been found, and we 
will not be
interested in the value of the $O$-constant. Informally, a 
problem has a good
algorithm if an instance of the problem is almost as 
easily {\it solved\/} as
it is {\it formulated}.

Sometimes we will refer to a {\it probabilistic\/} 
algorithm, which is an
algorithm that may use a random number generator for 
drawing random bits. One
formalization of this is a {\it nondeterministic\/} Turing 
machine (see
[\VEMB]). Unless we use the word {\it probabilistic}, we 
do {\it not\/} allow
the use of random number generators, and if we wish to 
emphasize this we
talk of {\it deterministic\/} algorithms. In the case of a 
probabilistic
algorithm, the running time and the output are not 
determined by the input
alone, but both have, for each fixed value of the input 
data, a {\it
distribution}. The {\it expected\/} running time of a 
probabilistic
algorithm is the expectation of the running time for a 
given input. Studying
the complexity of a probabilistic algorithm means finding 
an upper bound for
the expected running time as a function of the length of 
the input.
For a few convenient rules that can be used for this 
purpose we refer to
[\LEPO, \S 12]. A probabilistic algorithm is called {\it 
good\/} if
its expected running time is $(l+2)^{O(1)}$, where $l$ is 
the length of the
input.

{\it Parallel\/} algorithms have not yet played any role 
in algorithmic
number theory, and they will not be considered here.

Many results in this paper assert that ``there exists'' an 
algorithm with
certain properties. In all cases, such an algorithm can 
actually be
exhibited, at least in principle.

All $O$-constants are absolute and effectively computable 
unless indicated
otherwise.

\subheading{{\rm2.2.} Encoding data}
As stated above, the input and the output of an algorithm 
consist of
finite sequences of nonnegative integers. However, in the 
mathematical
practice of thinking and writing about algorithms one 
prefers to work with
mathematical concepts rather than with sequences of 
nonnegative integers that
encode them in some manner. Thus, one likes to say that 
the input of an
algorithm is given by an algebraic number field rather 
than by the sequence
of coefficients of a polynomial that defines the field; 
and it is both shorter
and clearer to say that one computes the kernel of a 
certain endomorphism of
a vector space than that one determines a matrix of which 
the columns express
a basis for that kernel in terms of a given basis of the 
vector space. To
justify such a concise mode of expression we have to agree 
on a way of
encoding entities such as number fields, vector spaces, 
and maps between them
by means of finite sequences of nonnegative integers. That 
is one of the
purposes of the remainder of this section. Sometimes there 
is one obvious way
to do the encoding, but often there are several, in which 
case the question
arises whether there is a good algorithm that passes from 
one encoding to
another. When there is, we will usually not distinguish 
between the encodings,
although for practical purposes they need not be equivalent.

We shall see that the subject of encoding mathematical 
entities suggests
several basic questions, but we will not pursue these 
systematically.
We shall not do much more than what will be needed in 
later sections.

\subheading{{\rm2.3.} Elementary arithmetic}
By $\ZZ$ we denote the ring of integers. Adding a sign bit 
we can clearly use
nonnegative integers to represent
{\it all\/} integers. The traditional algorithms for 
addition and subtraction
take time $O(l)$, where $l$ is the length of the input. 
The ordinary
algorithms for multiplication and division with remainder, 
as well as the
Euclidean algorithm for the computation of greatest common 
divisors, have
running time $O(l^2)$. With the help of more sophisticated 
methods this can
be improved to $l^{1+o(1)}$ for $l\to\infty$ (see [\KNUT]).
An operation that is {\it not\/} known to be doable by 
means of a good
algorithm is decomposing a positive integer into prime
numbers (see [\KNUT, \LETI, \LLAN]), but there is a good 
probabilistic
algorithm for the related problem of deciding whether a 
given integer is
prime [\ADHU]. No good algorithms are known for the 
problem of recognizing
squarefree numbers and the problem of finding the largest 
square dividing a
given positive integer, even when the word ``good'' is 
given a less formal
meaning (see [\LLMF, \S2]).

For some algorithms a prime number $p$ is part of the 
input. In such a case,
the prime is assumed to be encoded by itself rather than 
that, for example,
$n$ stands for the $n$\<th prime. Since we know no good 
deterministic algorithm
for recognizing primes, it is natural to ask what the 
algorithm does if
$p$ is not prime or at least not known to be prime. Some 
algorithms may
discover that $p$ is nonprime, either because a known 
property of primes is
contradicted in the course of the computations, or because 
the algorithm
spends more time than it should; such algorithms may be 
helpful as primality
tests. Other algorithms may even give a nontrivial factor 
of $p$, which may
make them applicable as integer factoring algorithms. For 
both types of
algorithms, one can ask what can be deduced if the 
algorithm does
appear to terminate successfully. Does this assist us in 
proving that $p$ is
prime? What do we know about the output when we do not 
assume that $p$ is
prime? An algorithm for which this question has not been 
answered
satisfactorily is Schoof's algorithm for counting the 
number of points on an
elliptic curve over a finite field [\SCEC].

{\it Rational numbers\/} can be represented as pairs of 
integers in an
obvious manner, and all field operations can be performed 
on them in
polynomial time.

Let $n$ be a positive integer. The elements of the ring 
$\ZZ/n\ZZ$ are
assumed to be encoded as nonnegative integers less than 
$n$. The ring
operations can be performed in polynomial time. An {\it 
ideal}
$I\subset\ZZ/n\ZZ$ can be encoded either by means of its 
index
$d=[\ZZ/n\ZZ:I]$, which completely determines it and which 
can be any divisor
of $n$, or by means of a finite sequence of elements that 
generates $I$, or
by means of a single generator. An element of $I$ can be 
represented either
as an element of $\ZZ/n\ZZ$ that is divisible by $d$, or 
as an explicit
$\ZZ/n\ZZ$-linear combination of the given generators of 
$I$, or as an
explicit multiple of a single given generator. Using the 
extended Euclidean
algorithm one easily sees that one can pass from any of 
these encodings of
ideals and their elements to any other in polynomial time 
and that one can
likewise test inclusion and equality of given ideals. In 
particular, one can
decide in polynomial time whether a given nonzero element 
of $\ZZ/n\ZZ$ is a
unit, if so find its inverse, and if not so find a 
nontrivial divisor of $n$.
Taking $n=p$ to be prime we conclude that we can perform 
all field operations
in $\FF_p=\ZZ/p\ZZ$ in polynomial time.

\subheading{{\rm2.4.} Linear algebra}
Let $F$ be a field, and suppose that one has agreed
upon an encoding of its elements, as is the case when $F$ 
is the field $\QQ$
of rational numbers or the field $\FF_p$ for some prime 
number $p$ (see 2.3).
Giving a finite-dimensional vector space over $F$ simply 
means giving
a nonnegative integer $n$, which is the dimension of the 
vector space.
This number $n$ is to be given in {\it unary}, i.e., as a 
sequence $1$,
$1$, \dots, $1$ of $n$ ones, so that the length of the 
encoding is at least
$n$. This is because almost any algorithm related to a 
vector space of
dimension $n$ takes time at least $n$. The elements of 
such a vector space
are encoded as sequences of $n$ elements of $F$.
Homomorphisms between vector spaces are encoded as 
matrices. A
subspace of a vector space can be encoded as a sequence of 
elements that
spans the subspace, or as a sequence of elements that 
forms a basis of the
subspace, or as the kernel of a homomorphism from the 
vector space to
another one. For all fields $F$ that we shall consider the 
traditional
algorithms from linear algebra, which are based on 
Gaussian elimination, are
polynomial-time: algorithms that pass back and forth 
between different
representations of subspaces, algorithms that decide 
inclusion and equality
of subspaces, that form sums and intersections of 
subspaces, algorithms that
construct quotient spaces, direct sums, and tensor 
products, algorithms
for computing determinants and characteristic polynomials 
of endomorphisms,
and algorithms that decide whether a given homomorphism is 
invertible and if
so construct its inverse. The proofs are straightforward, 
the main problem
being to find upper bounds for the sizes of the numbers 
that occur in the
computations, for example when $F=\QQ$.

If one applies any of these algorithms to $F=\ZZ/p\ZZ$ 
without knowing that
$p$ is prime, then one either finds a nontrivial divisor 
of $p$ because
some division by a nonzero element fails, or the algorithm 
performs
successfully as if $F$ were a field. In the latter case it 
is usually easy to
interpret the output of the algorithm in terms of free 
$\ZZ/p\ZZ$-modules
(see [\BULE]), thus avoiding the assumption that $p$ be 
prime.

\subheading{{\rm2.5.} Finitely generated abelian groups}
Specifying a finitely generated
abelian group is done by giving a sequence of nonnegative 
integers
$d_1$, $d_2$, \dots, $d_t$; the group is then 
$\bigoplus_{i=1}^t\ZZ/d_i\ZZ$,
which enables us to represent the elements of the group by 
means of
sequences of $t$ integers. In our applications the group 
is usually either
finite (all $d_i>0$) or free abelian (all $d_i=0$).
To make the $d_i$ unique one may require that $d_i$
divides $d_{i+1}$ for $1\le i<t$; this can be accomplished 
in polynomial
time.  One should not require the $d_i$ to be {\it prime 
powers}, since that
is, for all we know, algorithmically hard to achieve. 
Starting from this
description of finitely generated abelian groups, one can 
encode maps and
subgroups in various ways that are reminiscent of 2.4 and 
that are left to the
imagination of the reader. He may also formulate the 
analogues of the problems
mentioned in 2.4 for the current case and construct good 
algorithms for them
using Hermite and Smith reduction of integer matrices (see 
[\HMCB]). The main
difficulty is to keep the intermediate numbers small.

\subheading{{\rm2.6.} Basis reduction}
In many cases a finitely generated free abelian group $L$ 
is equipped with a
bilinear symmetric map $L\times L\to\RR$ that induces a 
Euclidean structure
on $L_\RR=L\otimes_\ZZ\RR$; here $\RR$ denotes the field 
of real numbers. For
example, this is the case if $L$ is a subgroup of $\ZZ^n$, 
with the ordinary
inner product. It is also the case if $L$ is a finitely 
generated subgroup of
the additive group of an algebraic number field $K$ (see 
2.9), the bilinear
symmetric map in this case being induced by 
$(x,x)=\sum_\sigma|\sigma x|^2$,
where $\sigma$ ranges over the field homomorphisms from 
$K$ to the field $\CC$
of complex numbers. In such cases it is often desirable to 
find a {\it reduced
basis\/} of $L$ over $\ZZ$, i.e., a basis of which the 
elements are
``short'' in a certain sense. If the symmetric matrix that 
defines the bilinear
map on a given basis of $L$ is known to a certain 
accuracy, then a reduced
basis can be found by means of a {\it reduction 
algorithm}. The complexity of
such an algorithm depends on the precise notion of 
``reduced basis'' that one
employs. In [\LLLF] one finds a good reduction algorithm 
that will suffice
for our purposes. See [\KANN] for further developments.

\subheading{{\rm2.7.}  Rings}
We use the convention that rings have unit elements, that
a subring has the same unit element, and that ring 
homomorphisms preserve the
unit element. The {\it characteristic\/} $\chr A$ of a 
ring $A$ is the
nonnegative integer that generates the kernel of the 
unique ring
homomorphism $\ZZ\to A$. The group of units of a ring $A$ 
is denoted by $A^*$.
All rings in this paper are supposed to be {\it 
commutative\/}.

Almost any ring that we need to encode in this paper has 
an additive group
that is either finitely generated or a finite-dimensional 
vector space over
$\QQ$; for exceptions, see 2.11.
Such a ring $A$ is encoded by giving its underlying 
abelian group as in
2.5 or 2.4 together with the multiplication map $A\otimes 
A\to A$.
It is straightforward to decide in polynomial time whether 
the multiplication
map satisfies the ring axioms.

{\it Ideals\/} are encoded as subgroups or, equivalently, 
as kernels of ring
homomorphisms. There are good algorithms for computing the 
sum, product, and
intersection of ideals, as well as the ideal $I:J=\{x\in 
A\colon\;xJ\subset
I\}$ for given $I$ and $J$, and the quotient ring of $A$ 
modulo a given ideal.

A {\it polynomial\/} over a ring is always supposed to be 
given by means of
a complete list of its coefficients, including the zero 
coefficients; thus
we do not work with sparse polynomials of a very high 
degree.

Most {\it finite rings\/} that have been encountered in 
algorithmic number
theory ``try to be fields'' in the sense that one is 
actually happy to find
a zero-divisor in the ring. This applies to the way they 
occur in \S 4
and also to the application of finite rings in primality 
testing [\LEGA,
\BOHU]. Nevertheless, it seems of interest to study finite 
rings from an
algorithmic point of view for their own sake. Testing 
whether a given finite
ring is local can be done by a good probabilistic 
algorithm, but finding the
localizations looks very difficult. Testing whether it is 
reduced or
a principal ideal ring also looks very difficult, but 
there may be a good
algorithm for deciding whether it is quasi-Frobenius. I do 
not know whether
isomorphism can be tested in polynomial time. Many 
difficulties are
already encountered for finite rings that are 
$\FF_p$-algebras for some
prime number $p$. Two finite \'etale $\FF_p$-algebras can 
be tested for
isomorphism in polynomial time (cf.\ [\BULE]), but there 
is no known
good deterministic algorithm for finding the isomorphism 
if it exists;
if they are fields, there is, but the proof depends on 
ring theory (see
[\LEFI]).

\subheading{{\rm2.8.}  Finite fields}
Let $p$ be a prime number, $n$ a positive integer,
and $q=p^n$. A finite field $\FF_q$ of cardinality $q$ is 
encoded as a ring,
as in 2.7. This comes down to specifying $p$, $n$, as well 
as a system of
$n^3$ elements $a_{ijk}$ of $\FF_p$ with the property that 
there is a basis
$e_1$, $e_2$, \dots, $e_n$ of $\FF_q$ over $\FF_p$ such 
that $e_ie_j=
\sum_ka_{ijk}e_k$ for all $i$, $j$. We refer to [\LEFI] 
for a description of
good algorithms for various fundamental problems: 
performing the field
operations in a given finite field, as well as 
exponentiation and the
application of automorphisms; finding all subfields of a 
given finite field
$\FF_q$, finding the irreducible polynomial of a given 
element of $\FF_q$
over a given subfield, finding a primitive element of 
$\FF_q$, i.e.,\ an
element $\alpha\in\FF_q$ with $\FF_q=\FF_p(\alpha)$, 
finding a normal
basis of $\FF_q$ over a given subfield, and finding all 
field homomorphisms
and isomorphisms from a given finite field to another. 
Most of these
algorithms rely heavily on linear algebra.

Given a positive integer $p$ and a system of $n^3$ 
elements $a_{ijk}$ of
$\ZZ/p\ZZ$, how does one decide whether they specify a 
field $\FF_q$ as above? This
is at least as hard as testing $p$ for primality, for 
which no good
deterministic
algorithm is known. However, this is the {\it only\/} 
obstruction: there is
a good algorithm that given $p$ and the $a_{ijk}$ {\it 
either\/} shows that
they do {\it not\/} define a field, {\it or\/} shows that 
if $p$ is prime they
do. Namely, one runs the algorithms mentioned above for 
finding a
primitive element $\alpha$ and its minimal polynomial $f$ 
over $\ZZ/p\ZZ$,
just as if one is working with a field, and one verifies 
that the map sending
$X$ to $\alpha$ induces an isomorphism from 
$(\ZZ/p\ZZ)[X]/(f)$ to the
structure that one is working with; if this is not true, 
or if anything
went wrong during the course of the algorithm, one does 
not have a field;
if it is, then as a final test one decides whether $f$ is 
irreducible over
$\ZZ/p\ZZ$, which for prime $p$ can be done by means of a 
good algorithm (see
[\AKLT, \LEFL] and the references given there).

There are also problems for which no good algorithm is 
known. One is the
problem of {\it constructing\/} $\FF_{p^n}$ for a given 
prime $p$ and a
given positive integer $n$, or, equivalently, constructing 
an irreducible
polynomial $f\in\FF_p[X]$ of degree $n$; here $n$ is 
supposed to be given in
unary (cf. 2.4). If one accepts the generalized Riemann 
hypothesis then there
is a good algorithm for doing this [\ADLE]. There is also 
a good {\it
probabilistic\/} algorithm for this problem, and a 
deterministic algorithm
that runs in $\sqrt p$ times polynomial time [\SHOA].

An important problem, which will come up several times in 
this paper, is
the problem of factoring a given polynomial $f$ in one 
variable over a given
finite field $\FF_{p^n}$. No good algorithm is known for 
this
problem, even when the generalized Riemann hypothesis is 
assumed. There
does exist a good probabilistic algorithm and a 
deterministic
algorithm that runs in $\sqrt p$ times polynomial time 
[\SHOB]; if $p$ is
fixed, or smaller than the degree of $f$, then the latter 
algorithm is good.
There also exists a good algorithm that, given 
$f\in\FF_{p^n}[X]$, determines
the {\it factorization type\/} of $f$, i.e., the number of 
irreducible
factors and their degrees and multiplicities. We refer to 
[\LEFL] for a
further discussion.

Algorithmic problems relating to the multiplicative group 
of finite fields,
such as the discrete logarithm problem, are generally very 
difficult, see
[\ODDL, \POFR, \LLAN, \GORD, \SCHI, \LOVO].

\subheading{{\rm2.9.} Number fields}
By a {\it number field\/} or an {\it algebraic
number field\/} we mean in this paper a field extension 
$K$ of finite degree
of the field $\QQ$ of rational numbers. For the basic 
theory of algebraic
number fields, see [\LANG, \WEIS, \CAFR].

An algebraic number field $K$ is encoded as its underlying 
$\QQ$-vector space
together with the multiplication map $K\otimes_\QQ K\to 
K$, as in 2.7; in
other words, giving $K$ amounts to giving a positive 
integer $n$ and a system
of $n^3$ rational numbers $a_{ijk}$ that describe the 
multiplication in $K$
on a vector space basis of $K$ over $\QQ$ (cf.\ 2.8 
above). As in [\LEFI,
\S 2], one shows that the field operations in a number 
field can be
performed in polynomial time. Using standard arguments 
from field theory one
shows that there are good algorithms for determining the 
irreducible
polynomial of a given element of $K$ over a given subfield 
and for finding
a primitive element of $K$, i.e., an element $\alpha\in K$ 
for which
$K=\QQ(\alpha)$. It follows that giving a number field is 
equivalent to
giving an irreducible polynomial $f\in\QQ[X]$ and letting 
the field be
$\QQ[X]/f\QQ[X]$.

Polynomials in one variable with coefficients in an 
algebraic number field
can be factored into irreducible factors in polynomial 
time. This is done
with the help of basis reduction, see [\LLLF, \LANA, 
\AKFA, \AKFB]. We note
two consequences.

First of all, from the argument given in 2.8 one sees that 
there is a good
algorithm for deciding whether a given system of $n^3$ 
rational numbers
defines a number field. Secondly, given {\it two\/} number 
fields
$K=\QQ(\alpha)$ and $K'$, one can decide whether or not 
they are isomorphic,
and if so, find all isomorphisms, in polynomial time. To 
do this, one factors
the irreducible polynomial $f$ of $\alpha$ over $\QQ$ into 
irreducible factors
in the ring $K'[X]$, and one observes that the {\it 
linear\/} factors are in
bijective correspondence with the field homomorphisms 
$K\to K'$; such a
field homomorphism is an isomorphism if and only if the 
two fields have the
same degree over $\QQ$.

With $K=K'$ we see from the above that one can also 
determine all
automorphisms of $K$, and composing them one can make a 
complete
multiplication table for the group $\Aut K$ of field 
automorphisms of $K$,
all in polynomial time.

In the proof of 3.5 we shall see that all maximal proper 
subfields of a given
number field of degree $n$ can be found in polynomial 
time. Finding
{\it all\/} subfields is asking too much, since the number 
of subfields is
not polynomially bounded. I do not know whether all {\it 
minimal\/} subfields
different from $\QQ$ can be found in polynomial time, nor 
whether their
number is $n^{O(1)}$. Intersections and composites of 
given subfields can be
found by means of linear algebra.

We stress that for our algorithms the number field $K$ is
considered to be {\it variable\/} rather than fixed, and 
that we wish our
running time estimates to be uniform in $K$.

\subheading{{\rm2.10.} Orders}
An {\it order\/} in a number field $K$ of degree $n$ is a
subring $A$ of $K$ of which the additive group is 
isomorphic to $\ZZ^n$.
Among all orders in $K$ there is a unique maximal one, 
which is called the
{\it ring of integers\/} of $K$ and denoted by $\Oo$. The 
orders in $K$ are
precisely the subrings of $\Oo$ of finite additive index. 
The {\it
discriminant} $\Delta_A$ of an order $A$ with $\ZZ$-basis 
$\omega_1$,
$\omega_2$, \dots, $\omega_n$ is the determinant of the 
matrix
$(\Tr(\omega_i\omega_j))_{i,j}$, where $\Tr\colon\ 
K\to\QQ$ is the
trace map. The discriminant of every order is a nonzero 
integer. The
discriminant of $\Oo$ is also called the discriminant of 
$K$ over $\QQ$ and
is simply denoted by $\Delta$.

There are several ways of encoding an order $A$ in a 
number field $K$. One is
by specifying $A$ as a ring as in 2.7, which amounts to 
giving $n$ and a
system of $n^3$ integers $a_{ijk}$; from 
$A\otimes_\ZZ\QQ\cong K$ it follows
that
the same data also encode $K$. Another is by specifying 
$K$ as well as a
sequence of elements of $K$ that generates $A$ as a ring, 
or as an abelian
group. We leave it to the reader to check that there are 
good algorithms
for transforming all these encodings into each other.

Given a number field $K$ one can construct an order in $K$ 
in polynomial time,
as follows. Let $n^3$ rational numbers $a_{ijk}$ be given 
that describe the
multiplication on a $\QQ$-basis $e_1=1$, $e_2$, \dots,
$e_n$ for $K$, and let $d$ be the least common multiple of 
the
denominators of the $a_{ijk}$. Then $A=\ZZ+\sum_{i=2}^n\ZZ 
de_i$ is an
order in $K$. In many cases one knows the irreducible 
polynomial $f$ of a
primitive element $\alpha$ of $K$ over $\QQ$. If 
$f\in\ZZ[X]$, then
one can take for $A$ the ``equation order'' $\ZZ[\alpha]$, 
which as a
ring is isomorphic to $\ZZ[X]/f\ZZ[X]$. If $f$ does not 
belong to $\ZZ[X]$,
then one can either replace $\alpha$ by $m\alpha$ for a 
suitable positive
integer $m$, or use a little known generalization of the 
equation order,
namely, the ring
$$A=\lf\{\beta\in 
K\colon\;\beta\cdot\sum_{i=0}^{n-1}\ZZ\alpha^i\subset
        \sum_{i=0}^{n-1}\ZZ\alpha^i\rt\}.$$
To find a $\ZZ$-basis for this ring, let $m$ be the least 
positive
integer for which the polynomial $g=mf=\sum_{i=0}^na_iX^i$
has coefficients $a_i$ in $\ZZ$ (with $a_n=m$); then
$$A=\ZZ+\sum_{i=1}^{n-1}\ZZ\cdot
        \(\sum_{j=0}^{i-1}a_{n-j}\alpha^{i-j}\).$$ 
These are exactly the rings $A$ for which $\Spec A$ is 
isomorphic to a
``horizontal'' prime divisor of the projective line over 
$\ZZ$.
Many results that are known for equation orders have 
direct analogues
for rings of this type; for example, the discriminant of 
$A$ equals the
discriminant of $g$. 

Applying basis reduction to a given order $A$ as in 2.6, 
one can find a 
$\ZZ$-basis for $A$ with the property that the integers 
$a_{ijk}$ that
express multiplication in this basis satisfy 
$a_{ijk}=|\Delta_A|^{O(n)}$.
This shows that $A$ can be encoded by means of data of 
length
$O(n^4(2+\log|\Delta_A|))$, and that there is a good 
algorithm for
transforming a given encoding into one satisfying this 
bound. From the
inequality $n\le2(\log|\Delta_A|)/\log3$, which is valid 
for all $A\ne\ZZ$,
one sees that the bound is $(2+\log|\Delta_A|)^{O(1)}$. It 
is often
convenient to assume that the given encoding of $A$ 
satisfies this bound, and
to estimate running times in terms of $|\Delta_A|$.

Let $A$ be an order in a number field $K$ of degree $n$. 
By a {\it fractional
ideal\/} of $A$ we mean a finitely generated nonzero 
$A$-submodule of $K$.
The additive group of a fractional ideal is isomorphic to 
$\ZZ^n$.
One can compute with fractional ideals as with ideals (see 
2.7).

\subheading{{\rm2.11.} Local fields}
A {\it local field\/} is a locally compact,
nondiscrete topological field. Such a field is 
topologically isomorphic to
the field $\RR$ of real numbers, or to the field $\CC$ of 
complex numbers, or,
for some prime number $p$, to a finite extension of the 
field $\QQ_p$ of
$p$-adic numbers, or, for some finite field $E$, to the 
field $E((t))$ of
formal Laurent series over $E$. A local field is 
uncountable, which implies
that we have to be satisfied with specifying its elements 
only to a certain
precision. The discussion below is limited to the case 
that the field is
non-archimedean, i.e.,\ not isomorphic to $\RR$ or $\CC$.

The complexity theory of local fields has not been 
developed as systematically
as one might expect on the basis of their importance in 
number theory (see
[\CASS]). The first thing to do is to
develop algorithms for factoring polynomials in one 
variable to a
given precision; see [\CHIA, \BULE] and \S 4 below. Here 
the incomplete
solution of the corresponding problem over finite fields 
(see 2.8) causes a
difficulty; we are forced to admit probabilistic 
algorithms, or to allow the
running time to be $\sqrt p$ times polynomial time, where 
$p$ denotes the
characteristic of the residue class field, or to avoid the 
need for {\it
completely\/} factoring polynomials. Once one can factor 
polynomials, it is
likely that satisfactory algorithms can be developed for 
the calculation of
ramification indices and residue class field degrees of 
finite extensions of
non-archimedean local fields. Some further problems are 
mentioned at the end
of \S3.
\heading3. Galois groups\endheading
In this section we are concerned with the following problem.

\dfn{Problem 3.1}
Given an algebraic number field $K$ and a nonzero 
polynomial $f\in K[X]$,
determine the Galois group $G$ of $f$ over $K$. Can this 
be done in
polynomial time? \enddfn

In the sequel we will always assume that the polynomial 
$f$ is squarefree.
This can be accomplished by means of a good algorithm, 
which replaces $f$ by
$f/\gcd(f,f')$. We denote the degree of $f$ by $n$.

We should specify how we want the algorithm to describe 
$G$. One
possibility is to require that the algorithm comes up with 
a complete
multiplication table of a finite group that is isomorphic 
to $G$, but this
has an important shortcoming. Namely, the group may be 
very large in
comparison to the length of the input, and it may not be 
possible to
write down such a complete multiplication table in 
polynomial time, let alone
calculate it. If we insist on a complete multiplication 
table, then
``polynomial time'' in Problem 3.1 should be taken to 
mean: polynomial time in
the combined lengths of the input {\it plus\/} output. 
Theorem 3.2 below
shows that Problem 3.1 does in this sense have a 
polynomial time solution.

If we are interested in more efficient algorithms, we 
should look for a
more concise way of describing $G$. For this, we view $G$ 
as a permutation
group of the zeroes of $f$ rather than as an abstract 
group. Numbering the
zeroes we see that $G$ may be regarded as a subgroup of 
the symmetric group
$S_n$ of order $n!$; this subgroup is determined only up 
to conjugacy due to
the arbitrary choice of the numbering of the zeroes. 
Instead of asking for a
multiplication table of $G$ we shall ask for a list of 
elements of $S_n$ that
generate $G$. Every subgroup of $S_n$ has a system of at 
most $n-1$ generators
(see [\MINE, Lemma 5.2]), and these can be specified using 
$O(n^2\log n)$ bits.
This is bounded by a polynomial function of the length of 
the input, since the
latter is at least $n$.

This formulation of the problem still leaves something to 
be desired;
namely, we do not ask how the numbering of the zeroes of 
$f$ is related to
other ways in which zeroes of $f$ may be specified: for 
example, as complex
numbers to a certain precision, for a suitable embedding 
$K\to\CC$, or
similarly as $p$-adic numbers for a suitable prime number 
$p$, or as elements
of an abstractly defined splitting field or of one of its 
subfields. However,
even without such a refined formulation the problem 
appears to be hard enough.

It should be remarked that a set of generators of a 
subgroup $G$ of $S_n$
can be used to answer, in polynomial time, several natural 
questions about
$G$. For example, one can determine its order; one can 
decide whether a
given element of $S_n$ belongs to $G$; one can, for a 
given prime $p$,
determine generators for a Sylow $p$-subgroup of $G$; one 
can
find a composition series for $G$ and name the isomorphism 
types of its
composition factors; in particular, one can decide whether 
$G$ is solvable.
For more examples, proofs, and references, see [\KALU]. It 
may be that some of
the ideas that underlie this theory, which depends on the 
classification of
finite simple groups, will play a role in a possible 
solution of Problem 3.1.

The following result, due to Landau [\LANA], expresses 
that the possibility
that $G$ is very large is the only obstruction to finding 
a good algorithm
for Problem 3.1.
\thm{Theorem 3.2}
There is a deterministic algorithm that given $K$ and $f$ 
as in\/ 
Problem {\rm3.1}
and a positive integer $b$ decides whether the Galois 
group $G$ has order at
most $b$, and if so gives a complete list of elements of 
$G$, and that runs
in time $(b+l)^{O(1)}$, where $l$ is the length of the 
data specifying $K$
and $f$.\ethm

The algorithm is obtained from the standard textbook 
construction
of a splitting field of $f$ over $K$. One first factors 
$f$ into
irreducible factors in $K[X]$. If all factors are linear, 
then the
splitting field is $K$ itself. Otherwise, one passes to 
the field
$L=K[X]/gK[X]$, where $g$ is one of the nonlinear 
irreducible factors of $f$.
Then a splitting field of $f$ over $L$ is also one over 
$K$, so applying
the algorithm recursively one can determine a splitting 
field of $f$ over $K$.
If at any stage during the recursion it happens that one 
obtains a
field that has degree larger than $b$ over the initial 
field $K$, then
$\#G>b$, and one stops. If this does not happen, then one 
eventually arrives
at a splitting field $M$ of $f$ over $K$. As in 2.9 one 
can determine
the group $\Gal(M/K)$ of all $K$-automorphisms of $M$, and 
this is $G$.
It is then easy to make a multiplication table for $G$ and 
to find an
embedding of $G$ into the symmetric group of the set of 
zeroes of $f$.

One sees from Theorem 3.2 that $G$ can be determined in 
time $(\#G+l)^{O(1)}$.
Since $\#G\le n!$, it follows that for bounded $n$ Problem 
3.1 is solved in
the sense that there is a polynomial time solution. This 
is an example of a
complexity result that does not adequately reflect the 
practical situation:
the practical problem of determining Galois groups is {\it 
not\/} considered
to be well solved, even though the algorithms that are 
actually used nowadays
always require $n$ to be bounded---in fact, each value of 
$n$ typically has
its own algorithm (cf.\ [\STAU, \FOMK]), which does {\it 
not\/} follow the
crude approach outlined above.
\proclaim {Corollary 3.3}
There is a good algorithm that given $K$ and $f$ decides 
whether $G$ is
abelian, and determines $G$ if $G$ is abelian and $f$ is 
irreducible.
\endproclaim

For irreducible $f$ this is easily deduced from Theorem 
3.2 with $b=n$, since a
transitive abelian permutation group of degree $n$ has 
order $n$. For
reducible $f$ one uses that the Galois group of $f$ is 
abelian if and only
if the Galois group of each irreducible factor of $f$ is 
abelian.

For reducible $f$, this algorithm does not determine the 
Galois group, and
it is not clear whether this can be done in polynomial 
time. The following
problem illustrates the difficulty.
\dfn{Problem 3.4}
Given an algebraic number field $K$ and elements $a_1$, 
$a_2$, \dots,
$a_t\in K$, determine the Galois group of\/ 
$\prod_{i=1}^t(X^2-a_i)$ over
$K$. Is there a good algorithm for doing this?\enddfn

For $K=\QQ$ this is indeed possible. For general algebraic 
number fields
one can probably do it if one assumes the generalized 
Riemann hypothesis.
Without such an assumption already the case that all $a_i$ 
are units of
the ring of integers of $K$ is difficult to handle. In any 
case, the
algorithm from Theorem 3.2 is in general too slow.

The following pretty result is due to Landau and Miller 
[\LAMI]. It shows
that one can decide in polynomial time whether $f$ is 
solvable by radicals
over $K$.
\proclaim{Corollary 3.5}
There is a good algorithm that given $K$ and $f$
decides whether $G$ is solvable.\endproclaim

As in the proof of Corollary
3.3, we may assume that $f$ is irreducible. If there were
a bound of the form $n^{O(1)}$ for the order of a solvable 
transitive
permutation group of degree $n$, then we could proceed in 
the same way as
for abelian groups. However, no such bound exists, since 
for every integer
$k\ge0$ there is a solvable transitive permutation group 
of degree $n=2^k$
and order $2^{n-1}$. Instead, one uses that the order of a 
{\it primitive\/}
solvable permutation group of degree $n$ does have an 
upper bound of the form
$n^{O(1)}$ (see [\PALF]). By Galois theory, the Galois 
group $G$ of $f$ is
primitive if and only if there are no nontrivial 
intermediate fields between
$K$ and $K(\alpha)$, where $f(\alpha)=0$. To reduce the 
general case to this
situation, it suffices to find a chain of fields 
$K=K_0\subset
K_1\subset\cdots\subset K_t=K(\alpha)$ that cannot be 
refined, since $G$ is
solvable if and only if for each $i$ the Galois closure of
$K_i\subset K_{i+1}$ has a solvable Galois group. Such a 
chain can be found
inductively if one can, among all intermediate fields 
$K\subset
L\subset K(\alpha)$ with $L\ne K(\alpha)$, find a maximal 
one. This is done
as follows. Factor the polynomial $f$ into monic 
irreducible factors over
$K(\alpha)$. One of the factors is $X-\alpha$. For each 
other irreducible
factor $g$ we define a subfield $L_g\ne K(\alpha)$ 
containing $K$
as follows. If $g$ is linear, $g=X-\beta$, then 
$K(\alpha)$ has a unique
$K$-automorphism $\sigma$ with $\sigma\alpha=\beta$, and 
we let $L_g$ be
the field of invariants of $\sigma$. If $g$ is nonlinear, 
then let $\beta$
be a zero of $g$ in an extension field of $K(\alpha)$, and 
 $L_g=
K(\alpha)\cap K(\beta)$. I claim that all maximal 
subfields are among the
$L_g$, so that we can find a maximal subfield by choosing 
a field $L_g$ with
the largest degree over $K$. The correctness of the claim 
follows by Galois
theory from the following purely group theoretic 
statement. Let $G$ be a
finite group, $H\subset J\subset G$ subgroups with $H\ne 
J$, and assume
that there is no subgroup $I$ of $G$ with $H\subset 
I\subset J$, $H\ne I\ne J$;
then there exists $\sigma\in G-H$ such that
$$\aligned
\langle H,\sigma\rangle=& J\quad\hbox{if}\ \sigma 
H\sigma^{-1}=H,\\
	\langle H,\sigma H\sigma^{-1}\rangle=& J\quad\hbox{if}
\ \sigma H\sigma^{-1}\ne H.
\endaligned$$
In fact, it suffices to choose $\sigma\in J-H$.

This concludes the sketch of the proof of Corollary 3.5. 
Note that the algorithm
does not determine the group $G$ if it is solvable, even 
if $f$ is
irreducible. One does obtain the prime divisors of $\#G$ 
if $G$ is
solvable.

Theorem 3.2 suggests that the largest groups are the 
hardest to determine.
However, the following result, which is taken from 
[\LANB], shows that the
{\it very\/} largest ones can actually be dealt with in 
polynomial time.
As above, let $S_n$ denote the full symmetric group of 
degree $n$, and let
$A_n$ be the alternating group of degree $n$.
\proclaim{Theorem 3.6}
There is a good algorithm that given $K$ and $f$ decides 
whether the
Galois group of $f$ is $S_n$ and whether or not it is 
$A_n$.\endproclaim

For this, one may by the above assume that $n\ge8$. From 
the classification
of finite simple groups it follows (see [\CAME]) that the 
only sixfold
transitive permutation groups of degree $n$ are $A_n$ and 
$S_n$. Hence, if we
build up the splitting field of $f$ over $K$ as in the 
proof of Theorem 3.2,
then $G$ is $A_n$ or $S_n$ if and only if after adjoining 
six zeroes of $f$
one has obtained an extension of degree 
$n(n-1)(n-2)(n-3)(n-4)(n-5)$. One
can distinguish between $A_n$ and $S_n$ by computing the 
discriminant
$\Delta_f$ of $f$---this comes down to evaluating a 
determinant, which can
be done in polynomial time---and checking whether 
$X^2-\Delta_f$ has a
zero in $K$.

In a similar way one can decide in polynomial time whether 
$G$ is doubly
transitive. If $G$ is doubly transitive, one can determine 
the isomorphism
type of the unique minimal normal subgroup of $G$ in 
polynomial time, a result
that is due to Kantor [\KANT]. If one attempts to 
determine $G$ itself, one
runs into the following problem, which was suggested by 
Kantor.
\dfn{Problem 3.7}
Is there a polynomial time algorithm that given $K$ and 
$f$ as in\/ 
Problem {\rm3.1}
and a prime number $p$ decides whether $G$ has a normal 
subgroup of index $p$?
\enddfn

Even for $p=2$ this appears to be difficult.

Resolvent polynomials, such as $X^2-\Delta_f$ in the proof 
of 
Theorem 3.6, play a
much more important role in practical algorithms for 
determining Galois
groups than in known complexity results (see [\STAU, 
\FOMK]).
\dfn{Problem 3.8}
Is there a way to exploit resolvent polynomials to obtain 
complexity
results for varying $n$?\enddfn

The results that we have treated so far are more algebraic 
than arithmetic
in nature, the only exception being what we said about 
Problem 3.4. It
should be possible to formulate and prove similar results 
for other
sufficiently explicitly given fields over which 
polynomials in one variable
can be factored efficiently. We now turn to techniques 
that do exploit
the arithmetic of the field. The natural way to do this is 
to first
consider the case of finite and local base fields.

Let $E$ be a finite field, $f\in E[X]$ a nonzero 
polynomial, and $n$
its degree. As we mentioned in 2.8, there is a good 
algorithm that, given $E$
and $f$, determines the factorization type of $f$ in 
$E[X]$. This immediately
gives rise to the Galois group $G$, which is cyclic of 
order equal to
the least common multiple of the degrees of the 
irreducible factors of $f$.
One also obtains the cycle pattern of a permutation that 
generates $G$ as
a permutation group. Note that already in the case of 
finite fields the order
of $G$ may, for reducible $f$, be so large that the 
elements of $G$ cannot
be listed one by one in polynomial time.

We next discuss local fields.
\dfn{Problem 3.9}
Given a local field $F$ and a polynomial $f\in F[X]$ with 
a nonzero
discriminant, determine the Galois group $G$ of $f$ over 
$F$. What is the
complexity of this problem? Is there a good algorithm for 
it?\enddfn

I am not aware of any published work that has been done on 
Problem 3.9, and
I will only make a few brief remarks, restricting myself 
to the case that
$F$ is non-archimedean. Once a satisfactory theory of 
factoring polynomials
has been developed (see 2.11), one can prove an analogue 
of Theorem 3.2.
This does not yet solve the problem, since even when $f$ 
is irreducible the
Galois group may have a very large order. Tamely ramified 
extensions are small,
however, which suggests that the following problem should 
be doable.
\dfn{Problem 3.10}
Given $F$ and $f$ as in\/ Problem {\rm 3.9}, with $F$ 
non-archimedean, decide
whether a splitting field of $f$ over $F$ is tamely 
ramified, and if so
determine its Galois group over $F$. Can this be done in 
polynomial time?
\enddfn

When this problem is solved, one is left with wildly 
ramified extensions,
which occur only if $p$ is small. In that case, one may 
first want to consider
the following problem, which looks harder than Problem 3.10.
\dfn{Problem 3.11}
Given $F$ and $f$ as in\/ Problem {\rm 3.9}, with $F$ 
non-archimedean, determine
the Galois group of the maximal tamely ramified 
subextension $M$ of a
splitting field of $f$ over $F$. Can this be done in 
polynomial time?
\enddfn

If $f$ is irreducible of degree $n$, then the field $M$ in 
Problem 3.11 has
degree at most $n^4$ over $F$. This follows from a 
group-theoretic 
argument that was shown to me by I.$\,$M. Isaacs.

Even when all local problems are completely solved it is 
not clear whether
they are very helpful in solving Problem 3.1. There is a 
well-known
heuristic technique that can be used to obtain information 
about the
Galois group, which comes down to first considering the 
local Galois group
at primes that do not divide the discriminant of $f$ (see 
[\VDLI, \S1]).
Not much can be proved about this method, however (cf.\ 
[\LANB, \S4]).
G. Cornell has suggested to look instead at the {\it 
ramifying\/} primes,
the rationale being that Problem 3.1 should be reducible 
to the case
$K=\QQ$, in which case the Galois group is generated by 
the inertia groups.
\heading 4. Rings of integers\endheading
In this section we consider the following problem and its 
complexity.
\dfn{Problem 4.1}
Given an algebraic number field $K$, determine its ring of 
integers $\Oo$.
\enddfn

Constructing an order in $K$ as in 2.10 we see that this 
problem is
equivalent to the following one.
\dfn{Problem 4.2}
Given an order $A$ in a number field $K$, determine the 
ring of integers
$\Oo$ of $K$.\enddfn

Much of the literature on this problem assumes that the 
given order is
an equation order $\ZZ[\alpha]$, and it is true that 
equation orders offer
a few advantages in the initial stages of several 
algorithms. It may be that
in many practical circumstances one never gets beyond 
these initial stages
(cf.\ [\BERW, Preface]), but in the worst case---which is 
what we are
concerned with when we estimate the complexity of a 
problem---these advantages
quickly disappear as the algorithm proceeds. For this 
reason we make no
special assumptions about $A$ except that it is an order.

Most of what we have to say about Problem 4.2 also applies 
to the following more
general problem.
\dfn{Problem 4.3}
Given a commutative ring $A$ of which the additive group 
is isomorphic to
$\ZZ^n$ for some $n$, and that has a nonvanishing 
discriminant over $\ZZ$,
determine the maximal order in $A\otimes_\ZZ\QQ$.\enddfn

It is not difficult to show that Problems 4.2 and 4.3 are 
equivalent under
deterministic polynomial time reductions.

The main result on Problem 4.1, which is due to Chistov 
[\CHIB,
 \BULE], is a negative
one.
\proclaim{Theorem 4.4}
Under deterministic polynomial time reductions, Problem\/ 
{\rm 4.1} is
equivalent to the problem of finding the largest square 
factor of a
given positive integer.\endproclaim

The problem of finding the largest square factor of a 
given positive integer
$m$ is easily reduced to Problem 4.1 by considering the 
number field $K=\QQ(\sqrt m)$.
For the opposite reduction, which in computer science 
language is a
``Turing'' reduction, we refer to the discussion following 
Theorem 4.6 below.

Since there is no known algorithm for finding the largest
square factor of a given integer $m$ that is significantly 
faster than
factoring $m$ (see [\LLMF, \S 2]), Theorem 4.4 shows that 
Problem 4.1 is
currently intractable. More seriously, even if someone 
{\it gives\/} us $\Oo$,
we are not able to recognize it in polynomial time, even 
if probabilistic
algorithms are allowed. Deciding whether the given order 
$A$ in Problem
4.2 equals
$\Oo$ is currently an infeasible problem, just as deciding 
whether a given
positive integer is squarefree is infeasible. This is not 
just true in
theory, it is also true in practice.

One possible conclusion is that $\Oo$ is not an object 
that one should
want to work with in algorithms. It may very well be that 
whenever $\Oo$
is needed one can just as well work with an order $A$ in 
$K$, and
{\it assume\/} that $A$ equals $\Oo$ until evidence to the 
contrary is
obtained. This may happen, for example, when a certain 
nonzero ideal of
$A$ is found not to be invertible; in that case one can, in
polynomial time, construct an order $A'$ in $K$
that strictly contains $A$ and proceed with $A'$ instead 
of $A$.

If it indeed turns out to be wise to avoid working with 
$\Oo$, then it is
desirable that more attention be given to general orders, 
both
algorithmically and theoretically (cf. [\SAND]). This is 
precisely
what has happened in the case of quadratic fields (cf.\ 
[\LERE, \LEPO, \HMCA]).

The order $A$ equals $\Oo$ if and only if all of its 
nonzero prime ideals
$\ppx$ are nonsingular; here we call $\ppx$ nonsingular if 
the local ring
$A_\ppx$ is a discrete valuation ring, which is equivalent 
to
$\dim_{A/\ppx}\ppx/\ppx^2=1$. One may wonder, if it is 
intractable to find
$\Oo$, can one at least find an order in $K$ containing 
$A$ of which the
singularities are bounded in some manner? One result of 
this sort is given
below in Theorem 4.7; it implies that given $A$, one can 
find an order $B$ in
$K$ containing $A$ such that all singularities $\ppx$ of 
$B$ are {\it plane\/}
singularities, i.e.,\ satisfy $\dim_{B/\ppx}\ppx/\ppx^2=2$.

The geometric terminology just used should remind us of a 
situation in which
there does exist a good method for finding the largest 
square factor, namely,
if we are dealing with polynomials in one variable over a 
field. Thus,
Theorem 4.4 suggests that, for a finite field $E$, finding 
the integral
closure of the polynomial ring $E[t]$ in a given finite 
extension of $E(t)$ is
a tractable problem, and results of this nature have 
indeed been obtained
(see [\CHIB]). In geometric language, this means that it 
is feasible to
resolve the singularities of a given irreducible algebraic 
curve over a given
finite field. The corresponding problem over fields of 
characteristic zero has
been considered as well (see [\TEIT]), and one may wonder 
whether the
geometric techniques that have been proposed can also be 
used in the context
of Problem 4.2. In any case, we can formulate Problem 4.2 
geometrically by asking for the
resolution of the singularities of a given irreducible 
{\it arithmetic curve}.

For many purposes, resolving singularities is a local 
problem, but as we
see from Theorem 4.4 that is not quite the case in the 
context of algorithms.
It may be that one only needs to look locally at those 
prime ideals
$\ppx$ of $A$ for which $\dim_{A/\ppx}\ppx/\ppx^2>1$, but 
how does one
 {\it find\/}
those prime ideals? And likewise, if 
$A\cong\ZZ[X]/f\ZZ[X]$ is an equation
order, then, as all textbooks point out, one only needs to 
look locally at
those prime numbers $p$ for which $p^2$ divides the 
discriminant of $f$, but
how does one {\it find\/} those prime numbers? By 
contrast, once one {\it knows\/}
at which $\ppx$ or $p$ to look, the problem does admit a 
solution. To
formulate it we introduce some notation.

Let $A$ be an order in a number field $K$ of degree $n$. Let
further $C$ be a subring of $A$; for us, the most 
interesting cases are
$C=A$ and $C=\ZZ$. For any nonzero prime ideal $\ppx$ of 
$C$ we define
$$A^{(\ppx)}=\{\beta\in\Oo\colon\;\ppx^m\beta\subset 
A\hbox{ for
some }m\in\ZZ_{\ge0}\};$$ this is the ``$\ppx$-primary 
part'' of $\Oo$ when
viewed modulo $A$. It is not difficult
to show that $A^{(\ppx)}$ is an order in $K$ and that it 
is the smallest
order in $K$ containing $A$ with the property that all its 
prime ideals
containing $\ppx$ are nonsingular. In addition, one has an 
isomorphism
$\Oo/A\cong\bigoplus_\ppx A^{(\germ p)}/A$ of $C$-modules, 
with $\ppx$ ranging over
the set of nonzero prime ideals of $C$, and $A^{(\ppx)}=A$ 
for all but
finitely many $\ppx$. Thus, to determine $\Oo$, it 
suffices to determine all
$A^{(\ppx)}$. For a single $\ppx$, we have the following 
result.
\proclaim {Theorem 4.5}
There is a good algorithm that given $K$, $A$, $C$,
$\ppx$ as above, determines~$A^{(\ppx)}$.\endproclaim

This is proved by analyzing an algorithm of Zassenhaus 
[\ZASS, \ZIMM]. We
briefly sketch the main idea. Let us first consider the 
case $C=\ZZ$. Denote
by $p$ the prime number for which $\ppx=p\ZZ$, and write 
$A^{(p)}=A^{(\ppx)}$.

One needs a criterion for $A$ to be equal to $A^{(p)}$. The
{\it multiplier ring} $R_\xaa$ of a nonzero $A$-ideal 
$\xaa$ is defined by
$$R_\xaa=\{\beta\in K\colon\;\beta\xaa\subset\xaa\};$$
this is an order in $K$ containing $A$.
By $\qqx$ we shall denote a typical prime ideal of $A$ 
that contains $p$,
and we let $\rrx$ be the product of all such $\qqx$. By 
standard commutative
algebra, $A$ equals $A^{(p)}$ if and only if all $\qqx$ 
are invertible, and
$\qqx$ is invertible if and only if $R_\qqx=A$. Also, each 
$R_\qqx$ is contained
in $R_\rrx$, so that we can decide whether or not $A$ 
equals $A^{(p)}$ by
looking at $R_\rrx$. More precisely, if $R_\rrx=A$ then 
$A=A^{(p)}$, and
if $R_\rrx$ properly contains $A$ then so does $A^{(p)}$, 
since clearly
$R_{\germ r}\subset A^{(p)}$.

I claim that to turn the above considerations into an 
algorithm it suffices
to have a way of determining $\rrx$. Namely, suppose that 
$\rrx$ is known.
Then one can determine $R_\rrx$ by doing linear algebra 
over $\FF_p$, using
that $pR_\rrx/pA$ is the kernel of the $\FF_p$-linear map
$A/pA\to\End(\rrx/p\rrx)$ that sends each $x\in A/pA$ to the
multiplication-by-$x$ map. If this map is found to be 
injective, then
$R_\rrx=A$, and the algorithm terminates with $A^{(p)}=A$. 
If it is not
injective, then $R_\rrx$ strictly contains $A$. In that 
case one replaces
$A$ by $R_\rrx$ and starts all over again. Note that the 
number of iterations
is bounded by $(\log|\Delta_A|)/(2\log p)$, where 
$\Delta_A$ denotes the
discriminant of $A$.

It remains to find an algorithm for determining $\rrx$. 
Since the ideals
$\qqx$ are pairwise coprime, $\rrx$ is their intersection, 
so $\rrx/pA$ is the
set of nilpotents of the finite ring $A/pA$. It can, again 
by linear algebra,
be found as the kernel of the $\FF_p$-linear map $A/pA\to 
A/pA$ that sends
each $x\in A/pA$ to $x^{p^t}$; here $t$ is the least 
positive integer for
which $p^t\ge n$.

This concludes the sketch of the algorithm underlying 
Theorem 4.5 for
$C=\ZZ$. For general $C$, one can either modify the above, 
or first
determine $A^{(p)}$ for $p=\chr C/\ppx$ and then find 
$A^{(\ppx)}$
inside $A^{(p)}$.

The above algorithm gives, with a few modifications, also 
something if $p$
is not supposed to be prime. This is expressed in the 
following theorem,
which is taken from [\BULE].
\proclaim {Theorem 4.6}
There is a good algorithm that given $K$ and $A$
as above, as well as an integer $q>1$, determines an order 
$B$ in $K$ that
contains $A^{(p)}$ for each prime number $p$ that divides 
$q$ exactly once.
\endproclaim

To prove this, one first observes that it suffices to 
exhibit a good
algorithm that given $K$, $A$ and $q$ either finds $B$ as 
in the statement of
the theorem, or finds a nontrivial factorization 
$q=q_1q_2$. Namely, in the
latter case one can proceed recursively with $q_1$ and 
$q_2$ to find orders
$B_1$, $B_2$, and one lets $B$ be the ring generated by 
$B_1$ and $B_2$.

To find $B$ or $q_1$, $q_2$, one applies the algorithm 
outlined above, with a
few changes. The first change is that one starts by 
checking that $q$ is not
divisible by any prime number $p\le n$; if it is, then 
either one finds a
nontrivial splitting of $q$, or $q$ is a small prime 
number and one can apply
the earlier algorithm. So let it now be assumed that $q$ 
has no prime factors
$p\le n$, and that $q>1$. The second change is that one 
replaces, in the above
algorithm, $p$ and $\FF_p$ everywhere by $q$ and 
$\ZZ/q\ZZ$. This affects the
linear algebra routines, which are only designed to work 
for vector spaces
over fields. However, as we indicated in 2.4, they work 
just as well for
modules over a ring $\ZZ/q\ZZ$, {\it until\/} some 
division in $\ZZ/q\ZZ$
fails, in which case one obtains a nontrivial factor $q_1$ 
of $q$. The third
change is that $\rrx/q\ZZ$ should now be calculated as the 
``radical of the
trace form,'' i.e.,\ as the kernel of the 
$\ZZ/q\ZZ$-linear map
$A/qA\to\Hom(A/qA,\ZZ/q\ZZ)$ that sends $x$ to the map 
sending $y$ to
$\Tr(xy)$, where $\Tr\colon\ A/qA\to\ZZ/q\ZZ$ is the trace 
map. If $q$ is a
prime number exceeding $n$ then this is the same $\rrx$ as 
above.

One can show that the modified algorithm has the desired 
properties, see
[\BULE]. This concludes our sketch of the proof of Theorem 
4.6.

Using Theorem 4.6 we can complete the proof of Theorem 
4.4. Namely, suppose
that one has an algorithm that determines the largest 
square divisor of any
given positive integer. Calling this algorithm a few 
times, one can determine
the largest squarefree number $q$ for which $q^2$ divides 
the discriminant of
$A$. Applying the algorithm of Theorem 4.6 to $q$ one 
obtains an order $B$
that contains $A^{(p)}$ for each prime $p$ for which $p^2$ 
divides the
discriminant of $A$, so that $B=\Oo$.

We now formulate a result that also gives information 
about the local
structure of $B$ at primes $p$ for which $p^2$ divides 
$q$. Let $A$ be an
order in a number field $K$, and let $q$ be a positive 
integer. We call $A$
{\it nonsingular\/} at $q$ if each prime ideal of $A$ 
containing $q$ is
nonsingular. We call $A$ {\it tame\/} at $q$ if for each 
prime ideal $\germ p$
of $A$ containing $q$ there exist an unramified extension 
$R$ of the ring
$\ZZ_p$ of $p$-adic integers, where $p=\chr A/\ppx$, a 
positive integer $e$
that is not divisible by $p$, and a unit $u\in R^*$, such 
that there is an
isomorphism
$$\lim_{\scriptstyle\leftarrow\atop\scriptstyle
	m}A/\ppx^m\cong R[X]/(X^e-uq)R[X]$$
of $\ZZ_p$-algebras. As a partial justification of the 
terminology, we remark
that for prime $q$ the order $A$ is tame at $q$ if and 
only if each prime
ideal $\germ p$ of $A$ containing $q$ is nonsingular and 
tamely ramified over
$q$; this follows from a well-known structure theorem for 
tamely ramified
extensions of $\ZZ_q$ (see [\WEIS, \S 3-4]).
If $A$ is tame at $q$ and $\ppx$ is a prime
ideal of $A$ containing $q$, then $\ppx$ is nonsingular if 
and only if either
$p=\chr A/\germ p$ divides $q$ exactly once or the number 
$e$ above equals $1$,
and otherwise $\germ p$ is a plane singularity.
\proclaim {Theorem 4.7}
There is a good algorithm that, given an order $A$
in a number field $K$ of degree $n$, finds an order $B$ in 
$K$ containing $A$
and a sequence of pairwise coprime divisors $q_i$, $1\le 
i\le t$, of the
discriminant of $B$, such that
\roster
\item"(i)"
$B$ is tame at $q=\prod_{i=1}^tq_i$\RM;
\item"(ii)"
all prime numbers dividing $q$ exceed $n$\RM;
\item"(iii)"
$B$ is nonsingular at all prime numbers $p$ that do not 
divide $q$.\endroster
\endproclaim

This follows from a closer analysis of the algorithm of 
Theorem 4.6. Using
this theorem and the properties of tameness, one can 
deduce the following
result, which expresses that one can approximate $\Oo$ as 
closely as can be
expected on the basis of Theorem 4.4.
\proclaim {Theorem 4.8}
There is a good algorithm that, given an order $A$
in a number field $K$, finds an order $B$ in $K$ 
containing $A$ and a positive
integer $q$ dividing the discriminant of $B$ such that 
$B=\Oo$ if and only if
$q$ is squarefree, and such that the primes dividing 
$[\Oo:B]$ are exactly
those that appear at least twice in $q$. Moreover, there 
is a good algorithm
that given this $B$ and a nontrivial square
dividing $q$ finds an order in $K$ that strictly 
contains~$B$.\endproclaim

Next we discuss an algorithm that does a little more than 
the algorithm
of Theorem 4.5. Namely, in addition to finding 
$A^{(\ppx)}$, it also finds all
prime ideals of $A^{(\ppx)}$ containing $\ppx$. It 
depends---not surprisingly,
if one considers the case of an equation order 
$\ZZ[\alpha]$---on an
algorithm for factoring polynomials in one variable over a 
finite field, see
2.8. Due to this ingredient it is not a deterministic 
polynomial time
algorithm any more, and it has no extension as Theorem 4.6 
that works for
nonprimes.
\proclaim {Theorem 4.9}
There is a probabilistic algorithm that runs in expected 
polynomial time, and
there is a deterministic algorithm that runs in 
$\sqrt{\chr C/\ppx}$
times polynomial time, that given $K$, $A$, $C$, $\ppx$ as 
in\/ Theorem
 {\rm4.5},
determine
\roster
\item"(i)" all prime ideals of $A$ containing $\ppx$\RM; 
\item"(ii)" the order $A^{(\ppx)}$\RM;
\item"(iii)" all prime ideals of $A^{(\ppx)}$ containing 
$\ppx$.\endroster\endproclaim

One can do part (i) by analyzing the structure of the 
finite ring $A/\germ p A$,
as the reader may check; below we give a different 
argument. Once one has (i),
one can do (ii) by Theorem 
4.5 and (iii) by applying (i) to $A^{(\ppx)}$. We sketch
an alternative way to proceed, in which one constructs 
$A^{(\ppx)}$ and the
prime ideals simultaneously without appealing to Theorem
4.5. Let it first be assumed
that $C=A$.

The algorithm works with a list of pairs $B$, $\qqx$ for 
which $B$ is an
order in $K$ with $A\subset B\subset A^{(\ppx)}$ and 
$\qqx$ is a prime ideal of
$B$ containing $\ppx$. Initially, there is only one pair 
on the list, namely,
$A$, $\ppx$. The purpose of the algorithm is to achieve 
that $\qqx$ is 
nonsingular as a prime ideal of $B$, for each pair $B$, 
$\qqx$ on the list. If
that happens, then $A^{(\ppx)}$ is the sum of all $B$'s, 
and, as it turns out,
the ideals $\qqx A^{(\ppx)}$ are pairwise distinct and are 
precisely all prime
ideals of $A^{(\ppx)}$ containing $\ppx$.

The algorithm deals with a given pair $B$, $\qqx$ in the 
following manner.
First one determines, by means of linear algebra over the 
finite field
$B/\qqx$, an element $\gamma\in K$ with $\gamma\notin B$, 
$\gamma\qqx\subset B$;
such an element exists, see [\WEIS, Lemma 4-4-3]. Next, 
one considers
$\gamma\qqx$. If $\gamma\qqx\not\subset\qqx$, then $\qqx$
 is nonsingular, and the
pair $B$, $\qqx$ is left alone. Suppose now that 
$\gamma\qqx\subset\qqx$. Then
$B[\gamma]$ is an order in $K$ in which $\qqx$ is an 
ideal, and using
linear algebra one determines the minimal polynomial $g$ of
$(\gamma\bmod\qqx)$ over the field $B/\qqx$. This 
polynomial is factored into
irreducible factors over $B/\qqx$. For each irreducible 
factor $(h\bmod\qqx)$
of $g$, one now adds the pair $B[\gamma]$, $\qqx+
h(\gamma)B[\gamma]$ to the
list, and one removes $B$, $\qqx$.

The above is repeated until all pairs are nonsingular.

If $C\ne A$, then one replaces the pair $C$, $\ppx$ by 
$A'=C+\ppx A$, $\ppx A$;
note that $\ppx A$ is a prime ideal of $A'$ with $A'/\ppx 
A\cong C/\ppx$.
Applying the above with $A'$ in the role of $A$ one finds 
the order
$A^{\prime(\ppx)}$ and all of its prime ideals containing 
$\ppx$. One easily shows
that $A^{(\ppx)}=A^{\prime(\ppx)}$, and intersecting the 
prime 
ideals just mentioned
with $A$ one finds (i). This concludes the sketch of the 
proof of Theorem 4.9.

We note that the above algorithm also gives a convenient 
way of evaluating
the valuations corresponding to the prime ideals 
containing $\ppx$. Namely, for
each nonsingular pair $B$, $\qqx$ the corresponding 
valuation $v$ is given by
$$v(\beta)=\max\{m\in\ZZ_{\ge0}\colon\;\gamma^m\beta\in 
B\}$$
for $\beta\in B$, $\beta\ne0$, where $\gamma$ is as 
constructed in the
algorithm.  Since each element of $K$ can be written as a 
quotient of
elements of $B$ this allows us to compute $v(\beta)$ for 
each $\beta\in K$.

It is well known that the $p$-adic valuations of a number 
field
$K=\QQ(\alpha)$ correspond bijectively to the irreducible 
factors of $f$
over $\QQ_p$, where $f$ is the irreducible polynomial of 
$\alpha$ over
$\QQ$. Thus Theorem 4.9 suggests that factoring 
polynomials in one variable over
$\QQ_p$ to a given precision can be done by a 
probabilistic algorithm that
runs in expected polynomial time and by a deterministic 
algorithm that
runs in $\sqrt p$ times polynomial time. A result of this 
nature is given
in [\BULE]; see also [\CHIA], where a more direct approach 
is taken.

We close this section with a problem that is geometrically 
inspired.
\dfn{Problem 4.10}
If all singularities of $A$ are plane singularities, can 
the algorithm
of\/ Theorem {\rm4.9} be arranged in such a way that the 
same applies to all rings
$B$ that are encountered?\enddfn

It may be of interest to see whether the methods that have 
been
proposed for the resolution of plane curve singularities 
[\BRKN, \TEIT] shed
any light on this problem. One may also wish to 
investigate the algorithm of
Theorem 4.6 from the same perspective.

An affirmative answer to Problem 
4.10 may improve the performance of the algorithm.
This is because the hypothesis on $A$ is often satisfied, 
for example, if $A$
is an equation order or a ``generalized'' equation order 
as in 2.10; and
finding $\gamma$ in the algorithm of Theorem 4.9 may 
become easier if $\qqx$
 is at
worst a plane singularity, so that it can be generated by 
two elements.
\heading 5. Class groups and units\endheading
In this section we discuss the following problem and its 
complexity.
\dfn {Problem 5.1}
Given an algebraic number field $K$, with ring of integers 
$\Oo$,
determine the unit group $\Oo^*$ and the class group 
$\Cl\Oo$ of $\Oo$.
\enddfn

First we make a few remarks on the statement of the 
problem. In the previous
section we saw that, given $K$, the ring $\Oo$ may be very 
hard to determine
and that consequently we may have to work with subrings 
$A$ of $\Oo$ that,
for all we know, may be different from $\Oo$. Thus, it 
would have been
natural to formulate the problem for any order $A$ in $K$ 
rather than just
for $\Oo$. We have not done so, for several reasons. The 
first is that only
very little work has been done for general orders in 
fields of degree greater
than $2$. The second is that most difficulties appear 
already in the case
 $A=\Oo$ and that some additional complications are 
avoided. Finally, it
is to be noted that all algorithms for calculating unit 
groups and class
groups that have been proposed are so time-consuming that 
the effort required
in determining $\Oo$ appears to be negligible in 
comparison; and it may very
well be that the best way of calculating the unit group 
and class group of
a general order $A$ proceeds by first determining $\Oo$, 
next calculating
$\Oo^*$ and $\Cl\Oo$, and finally going back to~$A$.

We shall denote by $n$ and $\Delta$ the degree and the 
discriminant of $K$
over $\QQ$. It will be assumed that $\Oo$ is given by 
means of a
multiplication table of length $(2+\log|\Delta|)^{O(1)}$, 
as in 2.10.
We shall bound the running times of the algorithms in 
terms of $|\Delta|$.

The next question to be discussed is how we wish $\Oo^*$ 
and $\Cl\Oo$ to be
specified. As an abstract group, we have $\Oo^*\cong
(\ZZ/w\ZZ)\oplus\ZZ^{r+s-1}$, where $w$ denotes the number 
of roots of
unity in $K$ and $r$, $s$ denote the number of real and 
complex archimedean
places of $K$, respectively. Determining $\Oo^*$ means 
specifying the images
of the standard generators of $(\ZZ/w\ZZ)\oplus\ZZ^{r+
s-1}$ under an
isomorphism to $\Oo^*$; and we also like to be provided 
with an algorithm
that calculates the inverse isomorphism. Using the 
logarithms at the infinite
places (see [\LANG, Chapter V, \S 1]) and basis reduction 
(see 2.6) one
can prove that both these things can be achieved if we 
have a set of
generators for $\Oo^*$. However, just {\it writing down\/} 
a set of generators
for $\Oo^*$ may be very time-consuming. Suppose, for 
example, that $K$ is real
quadratic, i.e.,\ $n=2$ and $\Delta>0$. Then $\Oo^*$ is 
generated by $-1$ and
a single unit $\varepsilon$ of infinite order. It is easy 
to see that the total
number of digits of the coefficients of $\varepsilon$ on 
the given basis
of $\Oo$ over $\ZZ$ equals $R(\log\Delta)^{O(1)}$, where 
$R$ denotes the
regulator of $K$; see [\LANG, Chapter V, \S 1] for the 
definition of
the regulator. It is reasonable to conjecture that, for an 
infinite
sequence of real quadratic fields, $R$ is as large as 
$\Delta^{1/2+o(1)}$.
Hence we cannot expect to be able to write down 
$\varepsilon$, let alone
calculate it, in time significantly less than 
$\Delta^{1/2}$. If we are
interested in more efficient algorithms, then units must 
be represented
in a different way, for example as a product 
$\prod\gamma^{k(\gamma)}$
of elements $\gamma\in K^*$ with integer exponents 
$k(\gamma)$ that may
be very large in absolute value. This leads to the 
question whether there
exists a system of generating units that one can express 
in this way using
substantially fewer than $|\Delta|^{1/2}$ bits. Also, the 
following problem
is suggested.
\dfn{ Problem 5.2}
Given a number field $K$, finitely many elements 
$\gamma\in K^*$, and, for
each $\gamma$, an integer $k(\gamma)\in\ZZ$, decide whether
$\varepsilon=\prod_\gamma\gamma^{k(\gamma)}$ is a unit, 
i.e., belongs to
$\Oo^*$, and whether it equals $1$. If it is a unit, then 
determine its
residue class modulo a given ideal and calculate, for a 
given embedding
$\sigma\colon\ K\to\CC$, the logarithm of 
$\sigma\varepsilon$ to a given precision.
\enddfn

It may be expected that the first of these---recognizing 
units---can be done
by means of a good algorithm, even when $\Oo$ is not given,
by means of {\it factor refinement\/} (cf.\ [\BASH]). Good 
results on the
other problems can probably be obtained with diophantine 
approximation
techniques, such as basis reduction (see 2.6).
The same applies to the following more general problem.
\dfn{Problem 5.3}
Given a number field $K$ and a finite set $\Gamma$ of 
elements $\gamma\in
K^*$, find sets of generators for the subgroups
$$
\biggl\{(k(\gamma))_{\gamma\in\Gamma}\in\ZZ^\Gamma\colon\;%
\prod_{\gamma\in\Gamma}
	\gamma^{k(\gamma)}=1\biggr\},\qquad
\biggl\{(k(\gamma))_{\gamma\in\Gamma}\in\ZZ^\Gamma\colon\;%
\prod_{\gamma\in\Gamma}
	\gamma^{k(\gamma)}\in\Oo^*\biggr\}
$$
of $\ZZ^\Gamma$ and calculate the regulator of the group 
of all units
of the form\/ 
$\prod_{\gamma\in\Gamma}\gamma^{k(\gamma)}$\!, 
$k(\gamma)\in\ZZ$,
to a given precision.\enddfn

Problems of this nature arise in several contexts:
 in an algorithm for factoring integers
[\LLMP, \BULP], in the discrete logarithm problem [\GORD, 
\SCHI], as we shall see below;
in the determination of unit groups and class groups.

Returning to Problem 5.1, we still have to describe how we 
wish the class group
$\Cl\Oo$ to be specified. It is a finite abelian group, so 
we may first of all
ask for positive integers $d_1$, $d_2$, \dots, $d_t$ such 
that there is an
isomorphism $\bigoplus_i\ZZ/d_i\ZZ\cong\Cl\Oo$ of abelian 
groups, and
secondly for ideals $\xaa_1$, $\xaa_2$, \dots, $\xaa_t$ 
such that one such
isomorphism sends the standard generators of 
$\bigoplus_i\ZZ/d_i\ZZ$ to
the ideal classes of the $\xaa_i$. Once the class group 
has been calculated in
this sense, it may remain very difficult to find the {\it 
inverse\/}
isomorphism: given an $\Oo$-ideal, to which ideal of the 
form
$\prod_i\xaa_i^{m(i)}$ is it equivalent? Even testing 
whether a given ideal is
{\it principal\/} may be very difficult.

The order $h=\#\Cl\Oo$ of the class group is bounded by
$|\Delta|^{1/2}(n+\log|\Delta|)^{n-1}$ (see Theorem 6.5). 
The example of
imaginary quadratic fields---i.e., $n=2$ and 
$\Delta<0$---shows that $h$
is often as large as
$|\Delta|^{1/2}(\log|\Delta|)^{O(1)}$. Hence, if we are 
willing to spend
time at least of order $|\Delta|^{1/2}$ then we could 
conceivably list all
ideal classes, and finding the inverse isomorphism might 
also become
doable.

The first thing to be discussed about Problem 
5.1 is whether it can be done at all,
efficiently or not. This is a question that is strangely 
overlooked
in most textbooks, two notable exceptions being [\BOSH] 
and [\CASS].
For the class group, one often finds the theorem that 
every ideal class
contains an integral ideal of norm at most the Minkowski 
constant
$(n!/n^n)(4/\pi)^s|\Delta|^{1/2}$, where $s$ denotes the 
number of
complex places of $K$. However, this does not show that 
the class group is
effectively computable if no effective procedure for 
deciding equivalence of
ideals is supplied.

We shall prove a theorem from which the effective 
computability of $\Oo^*$
and $\Cl\Oo$ is clear. We begin by introducing some 
notation. Let $K$ be a
number field of degree $n$ and discriminant $\Delta$ over 
$\QQ$. A {\it
place\/} $\ppx$ of $K$ is an equivalence class of 
nontrivial absolute values
of $K$. The set of archimedean places of $K$ is denoted by 
$S_\infty$. For
$\ppx\notin S_\infty$, the {\it norm\/} $\NNx\ppx$ of 
$\ppx$ is the cardinality
of the residue class field at $\ppx$. For each place 
$\ppx$, let
$|\;|_\ppx\colon\ K\to\RR_{\ge0}$ denote the unique 
absolute value belonging
to $\ppx$ with the property that $|2|_\ppx=2$ if $\ppx$ is 
real; $|2|_\ppx=4$ if
$\ppx$ is complex; and $|K^*|_\ppx=(\NNx\ppx)^\ZZ$ if 
$\ppx$ is non-archimedean.
The {\it height\/} $H(x)$ of an element
$x\in K$ is defined by 
$H(x)=\prod_\ppx\max\{1,|x|_\ppx\}$, the product
extending over all places $\ppx$ of $K$. For any set $S$ 
of places of $K$
with $S_\infty\subset S$ we let $K_S$ denote the group of 
$S$-units, i.e.,
the subgroup of $K^*$ consisting of those $x\in K^*$ that 
satisfy $|x|_\ppx=1$
for all places $\ppx$ of $K$ with $\ppx\notin S$; in 
particular, we have
$K_{S_\infty}=\Oo^*$ if $\Oo$ denotes the ring of integers 
of $K$.
\proclaim {Theorem 5.4}
Let $K$ be an algebraic number field, $\Delta$ its 
discriminant over $\QQ$,
and $s$ the number of complex places of $K$. Let
$d=(2/\pi)^s|\Delta|^{1/2}$, and 
$S=S_\infty\cup\{\ppx\colon\;\ppx$ is a
finite place of $K$ with $\NNx\ppx\le d\}$.
Then the group $K_S$ is generated by the set of those 
$x\in K_S$ for
which $H(x)\le d^2$, and the ideal class group of the ring 
of integers of
$K$ is generated by the ideal classes of the finite primes 
in $S$.\endproclaim

The proof of this theorem is given in \S 6.

\rem{Remark} The example of real quadratic fields shows 
that it is not
reasonable to expect that the group $K_{S_\infty}=\Oo^*$ 
is generated by
elements $x$ for which $H(x)$ is substantially smaller 
than $e^d$. The
group $K_S$ in Theorem 5.4 is generally much larger than 
$\Oo^*$, but it is
generated by elements that are much smaller.

The relevance of Theorem 5.4 for the effective 
determination of $\Oo^*$ and $\Cl\Oo$
comes from the exact sequence
$$0\to\Oo^*\to K_S\to\ZZ^{S-S_\infty}\to\Cl\Oo\to0.$$
The middle arrow sends an element $x\in K_S$ to the vector
$(\ord_\ppx x)_{\ppx\in S-S_\infty}$, where $\ord_\ppx x$ 
is the number of
factors $\ppx$ in $x$; so $|x|_\ppx=\NNx\ppx^{-\ord_\ppx 
x}$. The map
$\ZZ^{S-S_\infty}\to\Cl\Oo$ sends $(m(\ppx))_\ppx$ to the 
ideal class
of $\prod_\ppx\ppx^{m(\germ p)}$. The exactness at 
$\Cl\Oo$ follows from the
last assertion of Theorem 5.4, the exactness at the other 
places is clear.
\endrem

To calculate $\Oo^*$ and $\Cl\Oo$ from the sequence, one 
starts by
calculating the set of generators of $K_S$ given by 
Theorem 5.4. It is 
well known
that there are only finitely many elements of bounded 
height in $K$
(see [\SEMW, Chapter 2]), and from the proof of this 
result it is clear that
they can be effectively determined. Determining the prime 
ideal factorizations
of these generators one finds a matrix that describes the 
map
$K_S\to\ZZ^{S-S_\infty}$. Applying algorithms for finitely 
generated abelian
groups (see 2.5) one obtains $\Oo^*$ and $\Cl\Oo$ as the 
kernel and cokernel
of this map.

We now turn to complexity results for Problem 5.1. Most 
results that have been
obtained concern quadratic fields (see [\LERE, \SCLT, 
\HMCA]).
For general number fields, virtually all that is known can 
be found in
[\BUCA] (note that, in that paper, 
$R^{1/2}\Dd^\varepsilon$ in Theorem 2 is a
printing error for $R\Dd^\varepsilon$, and $\Dd^{1/2+
\varepsilon}$ in Theorem 4 is
a printing error for $R^{1/2}\Dd^\varepsilon$). The 
following theorem appears to
be true.
\proclaim{Theorem 5.5}
Given $K$ and $\Oo$, one can determine a set of generators 
of $\Oo^*$ and the
structure of $\Cl\Oo$ in time at most $(2+
\log|\Delta|)^{O(n)}|\Delta|^{3/4}$
by means of a deterministic algorithm and in expected time 
at most
$(2+\log|\Delta|)^{O(n)}|\Delta|^{1/2}$ by means of a 
probabilistic algorithm.
\endproclaim

In [\BUCA] one finds a weaker version of this result, in 
which $n$ is
kept fixed. The more precise result should follow by 
combining [\BUCA] with
results that appear in [\BUSH].

The algorithm underlying Theorem 5.5, for which we refer 
to [\BUCA] and the
references given there, is not the same as the method for 
effectively
determining $\Oo^*$ and $\Cl\Oo$ that we just indicated. 
However, there
does exist a connection between the two methods. Namely, 
the proof of
Theorem 5.4 depends on a lemma from combinatorial group 
theory that constructs
a set of generators of a subgroup $H$ of a group $G$ from 
a set of
generators of $G$ itself (see Lemma 6.3), 
whereas the algorithm of Theorem 5.5 constructs
generators of the group $\Oo^*$ by letting it act on a 
certain graph; and it
is well known that these two subjects are closely related 
(see [\SEAA]). It
would be of interest to understand this connection better, 
and to see
whether Theorem 5.5 can be deduced from a suitable version 
of Theorem 5.4.

The higher exponent $3/4$ in Theorem 5.5 in the case of a 
deterministic algorithm
is due to the use of algorithms for factoring polynomials 
over finite fields
(see 2.8). It suggests the following problem.
\proclaim {Problem 5.6}
\RM{Can the exponent $3/4$ in\/ Theorem 
5.5 be replaced by}\/ $1/2$\RM?\endproclaim

For quadratic fields the answer is affirmative. It is 
likely that the
method by which this is shown, which is not completely 
obvious, carries over
to general number fields.

We close this section with an imprecise description of a 
probabilistic
technique for the solution of Problem 5.1.

Let the notation be as introduced before Theorem 5.4, and 
let $S$ consist of
the archimedean primes of $K$ and the non-archimedean 
primes of norm up to
a certain bound $b$. One supposes that one has a method of 
drawing elements
of $K_S$ that are ``random'' in a certain sense. For 
example, the method
might consist of drawing elements $x$ of $K$ whose 
coordinates on the given
vector space basis of $K$ over $\QQ$ are uniformly 
distributed over a
certain set of rational numbers, such as the positive 
integers up to a
certain bound, and keeping only those $x$ that are found 
to belong to $K_S$.

To determine the class group and the units, one draws 
elements of $K_S$ until
one has the feeling that the subgroup $H$ that they 
generate is equal to all
of $K_S$. One may get this feeling if the number of 
elements that have been
drawn is well over $\#S$, which is the minimal number of 
generators of $K_S$
as an abelian group, and if it happened several times in 
succession that a
newly drawn element of $K_S$ was found to belong to the 
subgroup generated
by the elements drawn earlier; if Problem 
5.3 has a satisfactory solution then this
can be tested. Assuming that $H=K_S$ one can determine 
$\Oo^*$ and $\Cl\Oo$,
as above, as the kernel and cokernel of the map 
$\phi\colon\ 
H\to\ZZ^{S-S_\infty}$ that sends $x$ to $(\ord_\ppx 
x)_{\ppx\in S-S_\infty}$.

In general, one does not know that $H=K_S$, so that 
$\ker\phi$ and $\coker\phi$
can only be conjectured to be $\Oo^*$ and $\Cl\Oo$, 
respectively. One does
know that there is an exact sequence
$$0\to\ker\phi\to\Oo^*\to 
K_S/H\to\coker\phi\to\Cl\Oo\to(\Cl\Oo)/C_S\to0,$$
where $C_S$ is the subgroup of $\Cl\Oo$ generated by the 
ideal classes of the
finite primes in $S$. The sequence shows that $H$ has 
finite index in $K_S$
if and only if the conjectured class group $\coker\phi$ is 
finite and the
$\ZZ$-rank of the conjectured unit group $\ker\phi$ mod 
torsion is the same
as it is for the true unit group $\Oo^*$, namely 
$\#S_\infty-1$. If $H$ has
infinite index in $K_S$ one should of course continue 
drawing elements of
$K_S$.

The information that one has about the relation between 
the conjectured class
group $\coker\phi$ and the true class group $\Cl\Oo$ is 
particularly meagre:
one has a group homomorphism $\coker\phi\to\Cl\Oo$, but 
neither its
injectivity nor its surjectivity is known. It is 
surjective if and only if
the ideal classes of the finite primes in $S$ generate the 
class group, and
results of this nature are known only if the bound $b$ 
that defines $S$ is at
least $|\Delta|^{1/2}$ times a constant depending on $n$. 
However, a
significant improvement is possible if one makes an 
unproved assumption.
Namely, Bach [\BACH, Theorem 4]  showed that if the 
generalized Riemann
hypothesis holds, then $\Cl\Oo$ is generated by the ideal 
classes
of the prime ideals of norm at most $12(\log|\Delta|)^2$. 
Hence if we assume
the generalized Riemann hypothesis then the map 
$\coker\phi\to\Cl\Oo$
is surjective for values of $b$ that are much smaller than 
$|\Delta|^{1/2}$.
If the map is surjective, then the above exact sequence 
shows that
$$h'R'=hR\cdot[K_S:H],\tag 5.7$$
where $h=\#\Cl\Oo$ and $R=\reg\Oo^*$ are the true class 
number and regulator,
and $h'=\#\coker\phi$ and $R'=\reg\ker\phi$ the 
conjectured ones; here we
assume that $H$ contains all roots of unity in $K$, which 
can easily be
accomplished [\POZA, \S 5.4]. Now suppose that we are able 
to estimate
$hR$ up to a factor $2$, i.e., that we can compute a 
number $a$ with
$a/2<hR<a$; if one assumes the generalized Riemann 
hypothesis this can
probably be done by means of a good algorithm, as in 
[\BUWI].
Then we see from (5.7) that $h'R'$ also satisfies 
$a/2<h'R'<a$ if and only if
$H=K_S$, and if and only if one has both $\ker\phi=\Oo^*$ 
and
$\coker\phi=\Cl\Oo$.

The above indicates that on the assumption of the 
generalized Riemann
hypothesis it may be possible to find a much faster 
probabilistic algorithm
for determining $\Oo^*$ and $\Cl\Oo$ than the algorithm of 
Theorem 
5.5. This leads
to the following problem
\dfn{Problem 5.8}
Assuming the truth of the generalized Riemann hypothesis, 
find
a probabilistic algorithm for\/ Problem 
{\rm5.1} that, for fixed $n$, runs in expected
time $$\exp\bigl(O((\log|\Delta|)^{1/2}(\log\log|%
\Delta|)^{1/2})\bigr),$$
the $O$-constant depending on $n$.\enddfn

Of course, one also wants to know how the running time 
depends on $n$, and
which value can be taken for the $O$-constant. For 
imaginary quadratic fields
Problem 5.8 has been solved [\HMCA]. For a partial 
solution in the general case,
see [\BUCB].
\heading 6. Explicit bounds\endheading
In the present section we prove a few explicit bounds on 
units and class
numbers of algebraic number fields, including Theorem 5.4.
Several proofs in this section are most naturally formulated
in terms of ideles, as in [\CAFR, Chapter II]. To stress 
the elementary
character of the arguments I have chosen to use more 
classical language.

We denote by $K$ an algebraic number field of degree $n$ and
discriminant $\Delta$ over $\QQ$, and by $r$ and $s$ the 
number of real and
complex places of $K$, respectively. We embed $K$ in
$K_\RR=K\otimes_\QQ\RR$, which, as an $\RR$-algebra, is 
isomorphic to
$\RR^r\times\CC^s$. We choose such an isomorphism, so that 
each element
$a\in K_\RR$ has $r+s$ coordinates $a_i$, of which the 
first $r$ are
real and the last $s$ complex. We put $n_i=1$ for $1\le 
i\le r$ and
$n_i=2$ for $r+1\le i\le r+s$. The {\it norm\/} $N\colon\ 
K_\RR\to\RR$
is defined by $Na=\prod_{i=1}^{r+s}|a_i|^{n_i}$.

Identifying each copy of $\CC$ with $\RR^2$ by mapping $x+
yi$ to $(x+y,x-y)$
we obtain an identification of $K_\RR$ with the 
$n$-dimensional Euclidean
space $\RR^n$. It is well known that this identification 
makes the ring
of integers $\Oo$ of $K$ into a lattice of determinant 
$|\Delta|^{1/2}$
in $K_\RR$, and more generally every fractional 
$\Oo$-ideal $\xaa$ into a
lattice of determinant $\NNx\xaa\cdot|\Delta|^{1/2}$, 
where $\NNx$ denotes the
ideal norm. We shall write $$d=(2/\pi)^s|\Delta|^{1/2}.$$

Let $S$ be a set of places of $K$ with $S_\infty\subset 
S$. By $I_S$ we
denote the group of fractional $\Oo$-ideals generated by 
the finite primes
in $S$, and by $K_S$, as in \S 5, the group $\{a\in 
K^*\colon\;\Oo a\in
I_S\}$. Denote by $i_S\colon\ K_S\to K_\RR^*\times I_S$ 
the embedding defined
by $i_Sa=(a,\Oo a)$. We give $K_\RR^*\times I_S$ the 
product topology, where
$I_S$ is discrete. For any compact set $B\subset 
K_\RR^*\times I_S$ the set
$B\cap i_SK_S$ consists of elements of bounded height and 
is therefore finite.
Hence $i_SK_S$ is discrete. Also, $i_SK_S$ is clearly 
contained in the
subgroup $V_S$ of $K_\RR^*\times I_S$ consisting of those 
pairs $(a,\xaa)$ for
which $Na=\NNx\xaa$.
\proclaim{Theorem 6.1}
Let $K$ be an algebraic number field, and let $S$ be a set 
of places of $K$
containing $S_\infty$ and containing all finite places 
$\ppx$ with
$\NNx\ppx\le d$, with $d$ as above. Let $V_S$ be as above, 
and denote by $F_S$
the set of all elements $(b,\bbx)\in V_S$ for which 
$\bbx\subset\Oo$,
$\NNx\bbx\le d$, and $|b_i|\le d^{1/n}$ for $1\le i\le r+
s$. Then $F_S$ is a
compact subset of $V_S$ and $V_S=F_S\cdot 
i_SK_S$.\endproclaim
\par\noindent
{\it Proof.}
The compactness of $F_S$  follows easily from the 
definition of $F_S$ and the
fact that $V_S$ is closed in $K_\RR\times I_S$. To prove 
the last assertion,
let $(a,\xaa)\in V_S$. Then $a\cdot\xaa^{-1}$ is a lattice 
of determinant
$Na\cdot|\Delta|^{1/2}\cdot\NNx\xaa^{-1}=|\Delta|^{1/2}$ 
in $K_\RR$.
By Minkowski's lattice point theorem there exists a 
nonzero element
$b\in a\xaa^{-1}$ with all $|b_i|\le d^{1/n}$. From $\Oo 
b\subset a\xaa^{-1}$
it follows that $\Oo b=a\xaa^{-1}\bbx$ for some integral 
$\Oo$-ideal $\bbx$.
Comparing determinants
we see that $Nb=\NNx\bbx$, so $\NNx\bbx\le d$. This 
implies that $\bbx\in I_S$,
so we have $(b,\bbx)\in F_S$. If we write $b=ac$ then $c$ 
is a nonzero element
of $\xaa^{-1}$, so $c\in K^*$. Since we also have $\Oo 
c=\xaa^{-1}\bbx\in
I_S$, we even have $c\in K_S$, so $(a,\xaa)=(b,\bbx)\cdot 
i_Sc^{-1}$.
 This proves
Theorem 6.1.

It follows from Theorem 6.1 that $V_S/i_SK_S$ is {\it 
compact}, if $S$ is as in the
theorem. This allows one to deduce the Dirichlet unit 
theorem and the
finiteness of the class number. Namely, take for $S$ the 
set of {\it all\/}
places of $K$. From the exact sequence $0\to 
V_{S_\infty}\to V_S\to I_S\to0$
one obtains an exact sequence
$$0\to V_{S_\infty}/i_{S_\infty}\Oo^*\to 
V_S/i_SK_S\to\Cl\Oo\to0,$$
where $\Oo^*$ and $\Cl\Oo$ are as in \S5. The map to 
$\Cl\Oo$ is
continuous if the latter is given the discrete topology. 
Thus the compactness
of $V_S/i_SK_S$ implies that 
$V_{S_\infty}/i_{S_\infty}\Oo^*$ is compact, which
is essentially a restatement of the Dirichlet unit 
theorem, and that $\Cl\Oo$
is finite. In the same way one proves that $V_S/i_SK_S$ is 
compact for {\it
every\/} set $S$ of primes containing $S_\infty$, not just 
for those from 
Theorem 6.1.

{}From the exact sequence and Theorem 
6.1 we see that every element of $\Cl\Oo$ is the
ideal class of an integral ideal $\bbx$ of norm at most 
$d$. This implies the
last assertion of Theorem 5.4. It also follows that 
$d\ge1$. The other
assertion of Theorem 5.4 is a special case of the 
following theorem, in which the
height $H$ is as defined in \S5:
$$H(a)=\germ N(\Oo+\Oo a)^{-1}\cdot\prod_{i=1}^{r+
s}\max\{1,|a_i|\}^{n_i}.$$
\proclaim {Theorem 6.2}
Let $K$, $S$ be as in\/ Theorem  {\rm 6.1}, with $S$ 
finite. Write
$m_S=\max\{\NNx\ppx\colon\;\ppx
\in S-S_\infty\}$ if $S\ne S_\infty$ and $m_S=1$ if
$S=S_\infty$. Then the group $K_S$ is generated by the set 
of those $a\in K_S$
satisfying $H(a)\le dm_S$ and also by the set of those 
$a\in\Oo\cap K_S$
satisfying $H(a)\le d^2m_S$.\endproclaim

For the proof we need a lemma from combinatorial group 
theory, as well as a
topological analogue.
\proclaim {Lemma 6.3}
Let $G$ be a group, $P$ a set of generators for $G$, and 
$H$ a subgroup of
$G$. Let $F$ be a subset of $G$ such that $G=FH$. Then $H$ 
is generated
by its intersection with 
$F^{-1}PF=\{x^{-1}yz\colon\;x,z\in F,y\in P\}$.
\endproclaim
\par\noindent
{\it Proof}. Replacing $P$ by $P\cup P^{-1}$ we may assume 
that $P=P^{-1}$,
and replacing $F$ by a subset we may assume that the 
multiplication map
$F\times H\to G$ is bijective. Let $J\subset H$ be the 
subgroup generated by
$H\cap F^{-1}PF$. If $y\in P$, $z\in F$, then $yz=xh$ for 
some $x\in F$,
$h\in H$, and then $h=x^{-1}yz\in H\cap F^{-1}PF\subset 
J$. This proves that
$PF\subset FJ$, so $PFJ\subset FJ$. Hence the nonempty set 
$FJ$ is stable
under left multiplication by $P$, which by our assumptions 
on $P$ implies
that $FJ=G$. From $J\subset H$ and the bijectivity of 
$F\times H\to G$ we
now obtain $J=H$. This proves Lemma 6.3.
\proclaim {Lemma 6.4}
Let $G$ be a Hausdorff topological group, and denote by 
$G_1$ the connected
component of the unit element\/ $1$ of $G$. Let $P\subset 
G$ be a subset
containing\/ $1$ such that $G$ is generated by $P\cup G_1$.
Let $H\subset G$ be a discrete subgroup, and let $F$ be a 
compact subset of $G$
such that $G=FH$. Then $H$ is generated by its 
intersection with $F^{-1}PF$.
\endproclaim
\par\noindent
{\it Proof.} The set $H\cap F^{-1}F$ lies in the discrete 
subgroup $H$, so
$(G-H)\cup(H\cap F^{-1}F)$ is open, and it contains the 
compact set
$F^{-1}F$. Hence it contains $F^{-1}UF$ for some open 
neighborhood $U$ of
$1$. Intersecting with $H$ we see that $H\cap 
F^{-1}F=H\cap F^{-1}UF$.
The subgroup of $G$ generated by $U$ is open, so it 
contains $G_1$.
Therefore $G$ is generated by $P\cup U$. Applying Lemma 
6.3 we find that $H$ is
generated by
$$\eqalign{H\cap(F^{-1}(P{\cup}U)F)&=(H\cap 
F^{-1}PF)\cup(H\cap F^{-1}UF)\cr
&=(H\cap F^{-1}PF)\cup(H\cap F^{-1}F)=(H\cap 
F^{-1}PF),\cr}$$
where in the last step we use that $1\in P$. This proves 
Lemma 6.4.

To prove Theorem 6.2, we apply Lemma 6.4 to
$$\displaylines{G=V_S,\qquad H=i_SK_S,\qquad F=F_S,\cr
P=\{x\in V_S\colon\;x^2=1\}\cup\{((\NNx\ppx)^{1/n},\ppx)%
\colon\;\ppx\in
	S-S_\infty\},\cr}$$
where $F_S$ is as in Theorem  6.1 and where 
$(\NNx\ppx)^{1/n}$ 
is viewed as an element
of $K_\RR^*$ via the natural inclusion $\RR^*\subset 
K_\RR^*$. Using 
Theorem 6.1 one
readily verifies that the conditions of Lemma 6.4 are 
satisfied. Hence $K_S$ is
generated by the set of those elements $a\in K_S$ for 
which there exist
$(b,\bbx)$, $(c,\ccx)\in F_S$, and $(y,\yyx)\in P$ such that
$$(a,\Oo a)=(b,\bbx)^{-1}\cdot(y,\yyx)\cdot(c,\ccx).$$
Then $\Oo a=\bbx^{-1}\germ y\ccx$, so the denominator 
ideal $\den a$ of $a$
divides $\bbx$. For all $i\in\{1,2,\ldots,r+s\}$ we have
$$|b_i|\le d^{1/n},\qquad|y_i|\le m_S^{1/n},\qquad|c_i|\le 
d^{1/n},$$
so for each subset $J\subset\{1,2,\ldots,r+s\}$ we have
$$\gather
	\prod_{i\in J}|c_i|^{n_i}\le d^{n_J/n}
		\qquad\hbox{where }n_J=\sum_{i\in J}n_i,\\
	\prod_{i\in J}|b_i|^{-n_i}=
		\NNx\bbx^{-1}\cdot\prod_{i\notin J}|b_i|^{n_i}
		\le\NNx\bbx^{-1}\cdot d^{1-n_J/n},\\
	\prod_{i\in J}|a_i|^{n_i}=\prod_{i\in J}|b_i|^{-n_i}
		|c_i|^{n_i}|y_i|^{n_i}\le\NNx\bbx^{-1}\cdot d\cdot 
m_S.\endgather$$
Choosing $J=\{i\colon\;|a_i|>1\}$ we obtain
$$H(a)=\NNx(\den a)\cdot\prod_{i\in J}|a_i|^{n_i}\le
	\NNx\bbx\cdot\NNx\germ b^{x-1}\cdot d\cdot m_S=d\cdot 
m_S.$$
This proves the first assertion of Theorem 6.2. To prove 
the second assertion,
we use Minkowski's lattice point theorem to choose a 
nonzero element
$b'\in\germ b$ with $|b'_i|\le(d\cdot\NNx\bbx)^{1/n}$ for 
all $i$. Then
$b'\bbx^{-1}$ is an integral ideal of norm at most $d$, so 
$b'\in\Oo\cap K_S$.
Also $b'a\in\Oo\cap K_S$, and we have
$$\eqalign{H(b')&=\prod_i\max\{1,|b'_i|\}^{n_i}\le 
d\cdot\NNx\bbx\le d^2,\cr
	H(b'a)&\le\prod_i\max\{1,|b'_i|\}^{n_i}
			\cdot\prod_i\max\{1,|a_i|\}^{n_i}\cr
		&\le d\cdot\NNx\bbx\cdot\NNx\bbx^{-1}\cdot d\cdot 
m_S=d^2m_S.\cr}$$
Since we can write $a=(b'a)/b'$, this proves the second 
assertion of
Theorem  6.2.

\rem{Remark} Theorem 6.2 is also valid if the bound 
$d^2m_S$ is
 replaced by $\max\{d^2m_{S'},
dm_S\}$, where $S'=S_\infty\cup\{\ppx\colon\;\NNx\ppx\le 
d\}$. This is proved by
applying Theorem 6.2 to $S'$ and choosing a nonzero 
element of height at most
$d\cdot\NNx\ppx$ in each prime $\ppx\in S-S'$.
\endrem

As a further application of Theorem 6.1, we deduce upper 
bounds for the class
number $h=\#\Cl\Oo$ and for the product $hR$ of the class 
number and the
regulator $R=\reg\Oo^*$. The upper bound for $hR$ 
resembles the upper bound
that Siegel [\SIEG] proved using properties of the zeta 
function of $K$. For
similar upper bounds, see [\QUEM].
\proclaim {Theorem 6.5}
Let $K$ be an algebraic number field of degree $n$ and 
discriminant $\Delta$
over $\QQ$, and let $s$ denote the number of complex 
places of $K$. Let
$d=(2/\pi)^s|\Delta|^{1/2}$. Then the class number $h$ and 
the regulator $R$
of $K$ satisfy
$$\eqalign{h&\le d\cdot{(n-1+\log d)^{n-1}\over(n-1)!},\cr
	hR&\le d\cdot{(\log d)^{n-1-s}\cdot(n-1+\log 
d)^s\over(n-1)!}.\cr}$$
\endproclaim
\par\noindent
{\it Proof.} We saw above that every ideal class contains 
an integral ideal
of norm at most $d$, so
$$h\le\#\{\bbx\subset\Oo\colon\;\NNx\bbx\le d\}.$$
For each positive integer $m$, the number of $\Oo$-ideals 
of norm $m$ is at
most the number of vectors $x=(x_i)_{i=1}^n\in\ZZ_{>0}^n$ 
satisfying
$\prod_ix_i=m$. One proves this by considering how 
rational primes can split
in $K$. Thus we obtain
$$\#\{\bbx\subset\Oo\colon\;\NNx\bbx\le d\}\le\#
\lf\{x\in\ZZ_{>0}^n\colon\;\prod
	x_i\le d\rt\}.$$
Replacing each $x$ by the box $\prod_{i=1}^n(x_i-1,x_i]$ 
we can estimate the
right side by a volume:
$$\#\lf\{x\in\ZZ_{>0}^n\colon\;\prod x_i\le d\rt\}\le
	\vl\lf\{x\in\RR_{>0}^n\colon\;\prod\max\{1,x_i\}\le 
d\rt\}.$$
Writing $y_i=\log x_i$ we see that the volume is equal to 
$J(n,\log d)$,
where generally for $n\in\ZZ_{>0}$, $\delta\in\RR_{\ge0}$ 
we put
$$J(n,\delta)=\int_{y\in\RR^n,\sum_i\max\{0,y_i\}\le\delta}
		\exp\Bigl(\sum_iy_i\Bigr)\,dy.$$
This integral is found to be
$$\eqalign{J(n,\delta)&=e^\delta\cdot\sum_{i=0}^{n-1}{n-1%
\choose
		i}{\delta^i\over i!}\cr
	&\le e^\delta\cdot\sum_{i=0}^{n-1}{n-1\choose
                i}{(n-1)^{n-1-i}\delta^i\over(n-1)!}
		=e^\delta\cdot{(n-1+\delta)^{n-1}\over(n-1)!}.\cr}$$
Putting $\delta=\log d$ we obtain the inequality for $h$.

For $hR$, we apply Theorem  6.1 with $S$ equal to the set 
of all places of $K$.
Let $u=\#S_\infty-1=n-1-s$, and define the group 
homomorphism
$\lambda\colon\ V_S\to\RR^u\times I_S$ by
$\lambda(a,\xaa)=((n_i\log|a_i|)_{i=1}^u,\xaa)$. This is a 
surjective group
homomorphism with a compact kernel, so $\lambda i_SK_S$ is 
discrete in
$\RR^u\times I_S$ with a compact quotient. From the 
definition of the
regulator one derives that $hR$ equals the volume of a 
fundamental domain
for $\lambda i_SK_S$ in $\RR^u\times I_S$. Hence Theorem 
6.1 implies that
$hR\le\vl\lambda F_S$. For each nonzero $\Oo$-ideal $\bbx$ 
with
$\NNx\bbx\le d$ we have, by an easy computation,
$$\vl\lambda\{(b,\bbx)\in V_S\colon\;|b_i|\le 
d^{1/n}\hbox{ for all }i\}=
	{\bigl(\log(d/\NNx\bbx)\bigr)^u\over u!}.$$
Therefore
$$hR\le\sum_{\NNx\bbx\le 
d}{\bigl(\log(d/\NNx\bbx)\bigr)^u\over u!},$$
where the sum is over integral $\Oo$-ideals $\bbx$. 
Proceeding as with
$h$ one finds that this is bounded above by
$$\int_{y\in\RR^n,\sum_i\max\{0,y_i\}\le\delta}
	\exp\Bigl(\sum_iy_i\Bigr)\cdot
	{\bigl(\delta-\sum_i\max\{0,y_i\}\bigr)^u\over u!}\,dy,$$
with $\delta=\log d$. Using that $s=n-1-u\ge0$ one finds 
after some
computation the integral to be
$$e^\delta\cdot\sum_{i=0}^s{s\choose i}{\delta^{u+
i}\over(u+i)!}
	\le e^\delta\cdot\sum_{i=0}^s{s\choose i}
	{(n-1)^{s-i}\delta^{u+i}\over(n-1)!}=
	e^\delta\cdot{\delta^u\cdot(n-1+\delta)^s\over(n-1)!}.$$
This proves Theorem 6.5.

\rem{Remark} The upper bound for $h$ in Theorem 6.5 is 
also valid when $d$, at both
occurrences, is replaced by the Minkowski constant 
$d'=(n!/n^n)(4/\pi)^s
|\Delta|^{1/2}$ of $K$, since every ideal class contains 
an integral ideal
of norm at most $d'$.
 \endrem

\heading{Acknowledgments}\endheading

The author gratefully acknowledges the hospitality and 
support of the Institute
for Advanced Study (Princeton). I also thank Enrico
Bombieri, Johannes Buchmann, Joe Buhler, Gary Cornell, 
Pierre Deligne, Bas
Edixhoven, Boas Erez, David Ford, Marty Isaacs, Ravi 
Kannan, Bill Kantor, Susan
Landau, Andries Lenstra, Arjen Lenstra, Kevin McCurley, 
John McKay, Andrew
Odlyzko, Michael Pohst, Carl Pomerance, and Jeff Shallit 
for their assistance
and helpful advice.


\Refs

\ref\no{\ADHU}
\by L. M. Adleman and M. A. Huang 
\paper Recognizing primes in random polynomialtime  
\paperinfo Research report, Dept. of Computer Science, Univ. 
of Southern
California, 1988
\moreref \jour Lecture Notes in Math.
\publ Springer
\publaddr Heidelberg 
\toappear
\afterall 
Extended abstract: Proc. 19th Ann. ACM Sympos. on Theory 
of Computing (STOC),
ACM, New York 1987, pp. 462--469. \endref

\ref\no{\ADLE}
\by L. M. Adleman and H. W. Lenstra, Jr. 
\paper Finding irreducible polynomials over finite fields  
\inbook Proc. 18th Ann. ACM Sympos. on Theory of Computing 
(STOC)
\publ ACM
\publaddr New York
\nofrills\yr\nofrills (1986,
\pages 350--355 \endref

\ref\no{\ADPR}
\by L. M. Adleman, C. Pomerance, and R. S. Rumely 
\paper On distinguishing prime numbers from composite 
numbers  
\jour Ann. of Math. (2) 
\vol 117 
\yr 1983 
\pages 173--206
\endref

\ref\no{\ARCH}
\by Archimedes 
\paper The sand-reckoner  
\inbook in:  Opera quae quidem extant 
\publ J. Hervagius
\publaddr Basel
\yr 1544
\lang Greek and Latin\endref

\ref\no{\ATMO}
\by A. O. L. Atkin and F. Morain 
\book Elliptic curves and primality proving
\toappear \endref

\ref\no{\BACH}
\by E. Bach 
\paper Explicit bounds for primality testing and related 
problems 
\jour Math. Comp. 
\vol 55 
\yr 1990 
\pages 355--380 \endref

\ref\no{\BASH}
\by E. Bach and J. O. Shallit 
\paper Factor refinement  
\jour J. Algorithms 
\toappear
\endref

\ref\no{\BERW}
\by W. E. H. Berwick 
\book Integral bases  
\publ Cambridge Univ. Press
\publaddr Cambridge
\yr 1927 \endref

\ref\no{\BOSH}
\by Z. I. Borevi\v c and I. R. \v Safarevi\v c 
\book Teorija \v cisel 
\publ Izdat. ``Nauka''
\publaddr Moscow
\yr 1964
\moreref   \transl English transl.:  
\nofrills\book Number theory  
\publ Academic Press
\publaddr New York
\yr 1966 \endref

\ref\no{\BOHU}
\by W. Bosma and M. P. M. van der Hulst 
\paper Primality proving with cyclotomy  
\inbook Academisch proefschrift, Universiteit van Amsterdam
\yr 1990 \endref

\ref\no{\BRKN}
\by E. Brieskorn and H. Kn\"orrer 
\book Ebene algebraische Kurven  
\publ Birkh\"auser
\publaddr Basel
\yr 1981 \endref

\ref\no{\BUCA}
\by J. Buchmann 
\paper Complexity of algorithms in algebraic number theory 
\inbook Proceedings of the first
conference of the Canadian Number Theory Association
\ed R. A. Mollin
\publ  De Gruyter
\publaddr Berlin
\yr 1990
\pages 37--53 \endref

\ref\no{\BUCB}
\bysame \paper A subexponential algorithm for the 
determination of
class groups and regulators of algebraic number fields  
\inbook  S\'eminaire de Th\'eorie des Nombres, Paris 
1988--1989
\ed C. Goldstein
\publ Birkh\"auser
\publaddr Boston
\yr 1990
\pages 27--41 \endref

\ref\no{\BULE}
\by J. Buchmann and H. W. Lenstra, Jr.
\paperinfo Manuscript in preparation \endref

\ref\no{\BUSH}
\by J. Buchmann and V. Shoup \paper Constructing 
nonresidues in finite fields
and the extended Riemann hypothesis  
\paperinfo in preparation
\afterall Extended abstract:
Proc. 23rd Ann. ACM Sympos. on Theory of Computing (STOC), 
ACM, New York 1991, pp. 72--79. \endref

\ref\no{\BUWI}
\by J. Buchmann and H. C. Williams \paper On the 
computation of the class number
of an algebraic number field  
\jour Math. Comp. \vol 53
\yr 1989 
\pages\nofrills 679--688. \endref

\ref\no{\BULP}
\by J. P. Buhler, H. W. Lenstra, Jr., and C. Pomerance 
\paper Factoring integers
with the number field sieve  
\paperinfo in preparation \endref

\ref\no{\CAME}
\by P. J. Cameron \paper Finite permutation groups and 
finite simple
groups  
\jour Bull. London Math. Soc. 
\vol 13 
\yr 1981 
\pages 1--22 \endref

\ref\no{\CASS}
\by J. W. S. Cassels 
\book Local fields  
\publ Cambridge Univ. Press
\publaddr Cambridge
\yr 1986 \endref

\ref\no{\CAFR}
\by J. W. S. Cassels and A. Fr\"ohlich (eds.) 
\paper Algebraic number theory 
\inbook Proceedings of an instructional conference
\publ Academic Press
\publaddr London
\yr 1967 \endref

\ref\no{\CHIA}
\by A. L. Chistov 
\paper  Efficient factorization of polynomials over local 
fields  
\jour Dokl. Akad. Nauk SSSR 
\vol 293 
\yr 1987 
\pages 1073--1077
\moreref \transl English transl.: Soviet Math. Dokl. 
\vol 35 
\yr 1987 
\pages 430--433 \endref

\ref\no{\CHIB}
\bysame  \paper The complexity of constructing the ring of 
integers
of a global field  
\jour Dokl. Akad. Nauk SSSR 
\vol 306 
\yr 1989 
\pages 1063--1067
\moreref\nofrills English transl.: 
\jour Soviet Math. Dokl. 
\vol 39 
\yr 1989 
\pages 597--600 \endref

\ref\no{\COHE}
\by H. Cohen 
\book A course in computational algebraic number theory  
\bookinfo in preparation \endref

\ref\no{\COLA}
\by H. Cohen and A. K. Lenstra \paper Implementation of a 
new primality test 
\jour Math. Comp. 
\vol 48 
\yr 1987 
\pages 103--121 \endref

\ref\no{\COLH}
\by H. Cohen and H. W. Lenstra, Jr. \paper Primality 
testing and Jacobi sums 
\jour Math. Comp. 
\vol 42 
\yr 1984 
\pages 297--330 \endref

\ref\no{\FOMK}
\by D. J. Ford and J. McKay \paper Computation of Galois 
groups from
polynomials over the rationals   
\inbook Computer Algebra
\eds D. V. Chudnovsky and R. D. Jenks
\bookinfo Lecture Notes in Pure and Appl. Math.
\vol 113 
\publ Marcel Dekker
\publaddr New York
\yr 1989
\pages 145--150 \endref

\ref\no{\GORD}
\by D. Gordon 
\paper Discrete logarithms using the number field sieve
\toappear \endref %

\ref\no{\HMCA}
\by J. L. Hafner and K. S. McCurley \paper A rigorous 
subexponential algorithm
for computation of class groups  
\jour J. Amer. Math. Soc. 
\vol 2
\yr 1989 
\pages 837--850 \endref

\ref\no{\HMCB}
\bysame  \paper Asymptotically fast triangularization
of matrices over rings  
\jour SIAM J. Comput. 
\toappear \endref

\ref\no{\KANN}
\by R. Kannan \paper Algorithmic geometry of numbers  
\inbook Annual Review of Computer Sciences, vol. 2
\eds J. F. Traub, B. J. Grosz, B. W. Lampson, N. J. Nilsson
\publ Annual Reviews Inc.
\publaddr Palo Alto
\yr 1987
\pages 231--267 \endref

\ref\no{\KANT}
\by W. M. Kantor
\paperinfo unpublished \endref

\ref\no{\KALU}
\by W. M. Kantor and E. M. Luks \paper Computing in 
quotient groups 
\inbook Proc. 22nd Ann. ACM Sympos. on Theory of Computing 
(STOC)
\publ ACM
\publaddr\nofrills New York\yr 1990
\pages 524--534 \endref

\ref\no{\KNUT}
\by D. E. Knuth \paper The art of computer programming 
\paperinfo Vol. 2
\inbook Seminumerical Algorithms
\publ  Addison-Wesley
\publaddr Reading, MA
\yr second edition, 1981 \endref

\ref\no{\LANB}
\by S. Landau \paper Polynomial time algorithms for Galois 
groups 
\inbook Eurosam 84
\ed J. Fitch
\bookinfo Lecture Notes in Comput. Sci.
\vol 174
\publ Springer
\publaddr Berlin
\yr 1984
\pages 225--236 \endref

\ref\no{\LANA}
\bysame  \paper Factoring polynomials over algebraic 
number fields 
\jour SIAM J. Comput. 
\vol 14
\yr 1985 
\pages 184--195 \endref

\ref\no{\LAMI}
\by S. Landau and G. L. Miller \paper Solvability by 
radicals is in polynomial
time  
\jour J. Comput. System Sci. 
\vol 30 
\yr 1985
\pages 179--208 \endref

\ref\no{\LANG}
\by S. Lang 
\book Algebraic number theory  
\publ Addison-Wesley
\publaddr Reading, MA
\yr 1970 \endref

\ref\no{\AKLT}
\by A. K. Lenstra \paper Factorization of polynomials  
\inbook Computational Methods in Number Theory 
\eds H. W. Lenstra, Jr. and  R. Tijdeman
\bookinfo Math. Centre Tracts.
\vol 154/155
\publ Mathematisch Centrum Amsterdam
\yr 1982
\pages 169--198 \endref

\ref\no{\AKFA}
\bysame  \paper Factoring polynomials over algebraic 
number fields 
\inbook Computer Algebra 
\ed (J. A. van Hulzen, ed.)
\nofrills\bookinfo Lecture Notes in Comput. Sci.
\vol 162 
\publ Springer
\publaddr Berlin
\yr 1983
\pages 245--254 \endref

\ref\no{\AKFB}
\bysame  \paper Factoring multivariate polynomials over 
algebraic
number fields  
\jour SIAM J. Comput. 
\vol 16 
\yr 1987 \pages 591--598 \endref

\ref\no{\LLAN}
\by A. K. Lenstra and H. W. Lenstra, Jr. \paper Algorithms 
in number theory 
\inbook Handbook of Theoretical Computer Science,
Vol. A, Algorithms and Complexity
\ed J.~van Leeuwen
\publ Elsevier
\publaddr Amsterdam
\yr 1990
\pages 673--715 \endref

\ref\no{\LLLF}
\by A. K. Lenstra, H. W. Lenstra, Jr., and L. Lov\'asz
 \paper Factoring polynomials with rational coefficients  
\jour Math. Ann. \vol 261 
\yr 1982 
\pages 515--534 \endref

\ref\no{\LLMF}
\by A. K. Lenstra, H. W. Lenstra, Jr., M. S. Manasse, and 
J. M. Pollard
\book The factorization of the ninth Fermat number
\toappear \endref

\ref\no{\LLMP}
\bysame
 \paper The number field sieve  
\paperinfo in preparation
\afterall Extended abstract: Proc. 22nd
Ann. ACM Sympos. on Theory of Computing (STOC), 
ACM, New York 1990, pp. 564--572. \endref

\ref\no{\LERE}
\by H. W. Lenstra, Jr. \paper On the calculation of 
regulators and class
numbers of quadratic fields   
\inbook Journ\'ees Arithm\'etiques 1980
\ed J. Armitage
\bookinfo  London Math. Soc. Lecture Note Ser.
\vol 56  
\publ Cambridge Univ. Press
\publaddr Cambridge\yr 1982
\pages 123--150 \endref

\ref\no{\LEGA}
\bysame  \paper Galois theory and primality testing 
 \inbook  Orders and Their Applications
\eds I. Reiner, K. Roggenkamp
\bookinfo Lecture Notes in Math.
\vol 1142  
\publ Springer
\publaddr Heidelberg
\yr 1985
\pages 169--189 \endref

\ref\no{\LEFL}
\bysame  \paper Algorithms for finite fields  
 \inbook  Number Theory and Cryptography
\ed J. H. Loxton
\bookinfo London Math. Soc. Lecture Note Ser.
\vol 154  
\publ Cambridge Univ. Press
\publaddr Cambridge
\yr 1990
\pages 76--85 \endref

\ref\no{\LEFI}
\bysame  \paper Finding isomorphisms between finite fields  
\jour Math. Comp. 
\vol 56 
\yr 1991 
\pages 329--347 \endref

\ref\no{\LEPO}
\by H. W. Lenstra, Jr. and C. Pomerance 
\paper A rigorous time bound for factoring
integers  
\jour J. Amer. Math. Soc. 
\toappear \endref

\ref\no{\LETI}
\by H. W. Lenstra, Jr. and R. Tijdeman (eds.) 
\paper Computational methods in number theory  
\inbook Mathematical Centre Tracts
\vol 154/155  
\publ Mathematisch Centrum
\publaddr Amsterdam
 \yr 1982 \endref

\ref\no{\LOVO}
\by R. Lovorn \paper Rigorous, subexponential algorithms 
for discrete logarithms
over finite fields  
\paperinfo thesis, University of Georgia, in preparation 
\endref

\ref\no{\MINE}
\by A. McIver and P. M. Neumann 
\paper Enumerating finite groups  
\jour Quart. J. Math. Oxford Ser. (2) 
\vol 38 
\yr 1987 
\pages 473--488 \endref

\ref\no{\ODDL}
\by A. M. Odlyzko \paper Discrete logarithms in finite 
fields and their
cryptographic significance   
\inbook Advances in Cryptology
\eds T. Beth, N. Cot, and I. Ingemarsson
\bookinfo   Lecture Notes in Comput. Sci. \vol 209 
\publ Springer
\publaddr Berlin \yr 1985 \pages 224--314 \endref

\ref\no{\PALF}
\by P. P. P\'alfy \paper A polynomial bound for the orders 
of primitive solvable
groups  
\jour J. Algebra 
\vol 77 
 \yr 1982 
\pages 127--137 \endref

\ref\no{\POHS}
\by M. E. Pohst \paper Three principal tasks of 
computational algebraic number
theory   
\inbook Number Theory and Applications
\ed R. A. Mollin
\publ Kluwer
\publaddr Dordrecht \yr 1989 \pages 123--133 \endref

\ref\no{\POZA}
\by M. Pohst and H. Zassenhaus \book Algorithmic algebraic 
number theory 
\publ Cambridge Univ. Press
\publaddr Cambridge
\yr 1989 \endref

\ref\no{\POFR}
\by C. Pomerance \paper Fast, rigorous factorization and 
discrete logarithm
algorithms \inbook Discrete Algorithms and Complexity
\eds D. S. Johnson, T. Nishizeki, A. Nozaki, and H. S. Wilf
\publ Academic Press \publaddr Orlando
\yr 1987 \pages 119--143 \endref

\ref\no{\QUEM}
\by R. Qu\^eme \paper Relations d'in\'egalit\'es 
effectives en th\'eorie
alg\'ebrique des nombres  
\jour  S\'em. Th\'eor. Nombres Bordeaux,
\yr\nofrills 1987--1988,
\pages 19-01--19-19 \endref

\ref\no{\SAND}
\by J. W. Sands \paper Generalization of a theorem of 
Siegel  
\jour Acta Arith.
\vol 58 
\yr 1991  \pages 47--57 \endref

\ref\no{\SCHI}
\by O. Schirokauer 
\paper Discrete logarithms and local units  
\paperinfo in preparation \endref

\ref\no{\SCLT}
\by R. J. Schoof 
\paper Quadratic fields and factorization  
\inbook Computational Methods in Number Theory 
\eds H. W. Lenstra, Jr. and  R. Tijdeman
\bookinfo Math. Centre Tracts. 
\vol 154/155
\publ Mathematisch Centrum Amsterdam \yr 1982
\pages 235--286 \endref

\ref\no{\SCEC}
\bysame  \paper Elliptic curves over finite fields and the 
computation
of square roots \RM{mod}\,$p$  
\jour Math. Comp. 
\vol 44 
\yr 1985  \pages 483--494 \endref

\ref\no{\SEAA}
\by J-P. Serre 
\paper Arbres, amalgames, ${\roman{ SL}}_2$  
\jour Ast\'erisque \vol 46 
\yr 1977 \endref

\ref\no{\SEMW}
\bysame  \paper Lectures on the Mordell-Weil theorem  
\publ Vieweg
\publaddr Braunschweig
\finalinfo 1989 \endref

\ref\no{\SHAN}
\by D. Shanks \paper The infrastructure of a real 
quadratic field and its
applications  
\inbook Proceedings of the 1972 number theory conference,
Univ. 
of Colorado, Boulder, CO
\yr 1972 \pages 217--224 \endref

\ref\no{\SHOA}
\by V. Shoup \paper New algorithms for finding irreducible 
polynomials over
finite fields  
\jour Math. Comp. 
\vol 54  \yr 1990  \pages 435--447 \endref

\ref\no{\SHOB}
\bysame   \paper On the deterministic complexity of 
factoring polynomials
over finite fields  
\jour Inform. Process. Lett. \vol 33  
\yr 1990  \pages 261--267 \endref

\ref\no{\SIEG}
\by C. L. Siegel \paper Absch\"atzung von Einheiten  
\inbook Nachr. Akad. Wiss. G\"ottingen Math.-Phys. Kl. 
\nofrills\yr 1969 \pages 71--86
\moreref\ 
\paper Gesammelte Abhandlungen
\inbook  Band IV
\publ Springer
\publaddr Berlin
\yr 1979
\pages 66--81 \endref

\ref\no{\STAU}
\by R. P. Stauduhar \paper The determination of Galois 
groups  
\jour Math. Comp.
\vol 27  \yr 1973  \pages 981--996 \endref

\ref\no{\SZPI}
\by L. Szpiro \paper Pr\'esentation de la th\'eorie 
d'Arak\'elov 
\inbook  Current Trends in Arithmetical Algebraic Geometry 
\ed K. A. Ribet
\bookinfo Contemp. Math. \vol 67 
\publ Amer. Math. Soc.
\publaddr Providence, RI
\yr 1987 \pages 279--293 \endref

\ref\no{\TEIT}
\by J. Teitelbaum 
\paper The computational complexity of the resolution of 
plane
curve singularities  
\jour Math. Comp. \vol 54 \yr 1990  \pages 797--837 \endref

\ref\no{\TZDW}
\by N. Tzanakis and B. M. M. de Weger \paper How to solve 
explicitly a
Thue-Mahler equation \toappear \endref

\ref\no{\VDLI}
\by F. J. van der Linden \paper The computation of Galois 
groups 
\inbook Computational Methods in Number Theory             
\eds H. W. Lenstra, Jr. and  R. Tijdeman
\bookinfo Math. Centre Tracts. \vol 154/155
\publ Mathematisch Centrum Amsterdam
\yr 1982 \pages 199--211 \endref

\ref\no{\VEMB}
\by P. van Emde Boas \paper Machine models, computational 
complexity and number
theory  
\inbook Computational Methods in Number Theory 
\eds H. W. Lenstra, Jr. and  R. Tijdeman
\bookinfo Math. Centre Tracts. \vol 154/155
\publ Mathematisch Centrum Amsterdam
\yr 1982 \pages 7--42 \endref

\ref\no{\WEIS}
\by E. Weiss \book Algebraic number theory  
\publ McGraw-Hill
\publaddr New York \yr 1963 \moreref
reprinted
\publ Chelsea
\publaddr New York
\yr 1976 \endref

\ref\no{\ZANT}
\by H. Zantema 
\paper Class numbers and units  
\inbook Computational Methods in Number Theory 
\eds H. W. Lenstra, Jr. and  R. Tijdeman
\bookinfo Math. Centre Tracts. \vol 154/155
\publ Mathematisch Centrum Amsterdam
\yr 1982 \pages 212--234 \endref

\ref\no{\ZASS}
\by H. Zassenhaus 
\paper Ein Algorithmus zur Berechnung einer Minimalbasis 
\"uber
gegebener Ordnung  
\inbook Funktionalanalysis, Approximationstheorie, 
numerische 
Mathematik, Oberwolfach
1965
\eds L. Collatz, G. Meinardus, and H. Unger
\publ Birkh\"auser
\publaddr Basel
 \yr 1967 \pages 90--103 \endref

\ref\no{\ZIMM}
\by H. G. Zimmer \paper Computational problems, methods 
and results in
algebraic number theory  
\inbook Lecture Notes in Math. \vol 262  
\publ Springer
\publaddr Berlin \yr 1972 \endref

\endRefs
\enddocument